\documentclass[12pt]{article}
\usepackage{textcomp}
\usepackage{graphicx}           
\usepackage{epstopdf}
\usepackage{color}              
\usepackage{xcolor,graphicx}
\usepackage{geometry}
\usepackage{float}
\usepackage{times}
\usepackage{indentfirst}        
\usepackage{amsmath,amssymb,bm} 
\usepackage{amsthm}
\usepackage{cases}              
\usepackage{setspace}
\usepackage{hyperref}
\usepackage{titlesec}
\usepackage[all]{xy}            
\usepackage{tikz}
\usepackage{multirow}
\usepackage{graphicx} 
\usepackage{epstopdf}
\usepackage{caption}
\usepackage{diagbox}
\usepackage{mathrsfs}
\usepackage{tikz}
\usepackage{dsfont}
\usepackage{graphics}
\usepackage{hyperref,cleveref,cite}
\usepackage{tcolorbox}
\usepackage{fancyhdr,fancyvrb}
\usepackage{xcolor,graphicx,geometry,float,indentfirst,cases,setspace,titlesec}
\usepackage{diagbox}

\newtheorem{thm}{Theorem}[section]
\newtheorem{defi}[thm]{Definition}

\begingroup
\makeatletter \@addtoreset{equation}{section} \makeatother
\makeindex 
\def\qed{\hfill \rule{4pt}{7pt}}
\def\pf{\vskip 0.2cm {\noindent \bf Proof.}\quad}

\begin{document}

\begin{center}
{\Large \bf  The $(\alpha,\beta)$-Eulerian Polynomials and Descent-Stirling\\[5pt] Statistics on Permutations}

\vskip 4mm
Kathy Q. Ji
\vskip 2mm

Center for Applied Mathematics\\[2pt]
Tianjin University\\[2pt]
Tianjin 300072, P.R. China

\vskip 2mm
kathyji@tju.edu.cn

\end{center}

\vskip 6mm \noindent {\bf Abstract.} Carlitz and Scoville introduced the polynomials  $A_n(x,y|{\alpha},{\beta})$, which we refer to as the $(\alpha, \beta)$-Eulerian polynomials. These polynomials count permutations based on Eulerian-Stirling statistics, including descents, ascents, left-to-right maxima, and right-to-left maxima.   Carlitz and Scoville obtained the generating function of $A_n(x,y|{\alpha},{\beta})$.   In this paper, we introduce a new family of polynomials, $P_n(u_1,u_2,u_3,u_4|{\alpha},{\beta})$, defined on permutations, incorporating descent-Stirling statistics including valleys, exterior peaks, right double descents, left double ascents, left-to-right maxima, and right-to-left maxima.
 By employing the grammatical calculus introduced by Chen,  we establish the connection between the generating function of  $P_n(u_1,u_2,u_3,u_4|{\alpha},{\beta})$
and the generating function of the $(\alpha,\beta)$-Eulerian polynomials $A_n(x,y|{\alpha},{\beta})$ introduced by Carlitz and Scoville. Using this connection,  we derive the generating function of $P_n(u_1,u_2,u_3,u_4|{\alpha},{\beta})$, which can be specialized to obtain the $(\alpha,\beta)$-extensions of  generating functions for peaks, left peaks, double  ascents, right double  ascents and left-right double ascents given by David-Barton, Elizalde and Noy, Entringer, Gessel, Kitaev and Zhuang. Moreover,   we establish two   relations between $P_n(u_1,u_2,u_3,u_4|{\alpha},{\beta})$ and $A_n(x,y|{\alpha},{\beta})$, which enable us to derive $(\alpha,\beta)$-extensions of results of Stembridge, Petersen, Br\"and\'en, and Zhuang. We also obtain the left peak version of Stembridge's formula and the peak version of Petersen's formula, along with their respective $(\alpha,\beta)$-extensions, by utilizing these two relations.  Specializing $(\alpha,\beta)$-extensions of  Stembridge's formula and the left peak version of Stembridge's formula allows us to derive the $(\alpha,\beta)$-extensions of the tangent and secant numbers.



\section{Introduction}
The objective of this paper is to   investigate the  polynomials involving descent-Stirling statistics. Let us first recall Eulerian, Stirling and descent statistics on permutations.  Let $\mathfrak{S}_n$ denote the set of permutations on $[n]:=\{1,\,2,\ldots, n\}.$   We say that  $i$ is   a descent   of $\sigma=\sigma_1\sigma_2\cdots \sigma_n \in \mathfrak{S}_n$ if $1\leq i<n$ and $\sigma_{i}>\sigma_{i+1}$. The  Eulerian polynomials are defined by
\begin{equation}\label{def-Euler}
A_n(x):=\sum_{\sigma \in \mathfrak{S}_{n}} x^{{\rm des}(\sigma)+1}
\end{equation}
with the convention that $A_0(x)=1$, where ${\rm des}(\sigma)$ counts the number of descents of $\sigma$. The generating function of $A_n(t)$ is well known:
 \begin{equation}\label{gf-Euler}
 \sum_{n\geq 0} A_n(x) \frac{t^n}{n!}=\frac{1-x}{1-xe^{{(1-x)}t}}.
 \end{equation}
Eulerian polynomials carry a profound historical legacy and play a pivotal role across diverse combinatorial landscapes. For an extensive analysis, please refer to Petersen \cite{Petersen-2015}.

A permutation statistic whose generating function is given by \eqref{def-Euler} is called Eulerian.  Let $\sigma=\sigma_1\cdots \sigma_n \in \mathfrak{S}_n$. A number $1\leq i<n$ for which $\sigma_{i}<\sigma_{i+1}$ is called an ascent   of $\sigma$  and a number $1\leq i<n$ for which $\sigma_{i}>i$ is called an excedance   of $\sigma$. Let ${\rm asc}(\sigma)$ denote the number of ascents of $\sigma$  and let ${\rm exc}(\sigma)$ denote the number of excedances of $\sigma$. It is known that ${\rm asc}(\sigma)$ and ${\rm exc}(\sigma)$ are  Eulerian statistics, see MacMahon \cite[p.186]{MacMahon-1960} and  Stanley \cite[p.186]{MacMahon-1960}. To wit,
 \[A_n(x):=\sum_{\sigma \in \mathfrak{S}_{n}} x^{{\rm des}(\sigma)+1}=\sum_{\sigma \in \mathfrak{S}_{n}} x^{{\rm asc}(\sigma)+1}=\sum_{\sigma \in \mathfrak{S}_{n}} x^{{\rm exc}(\sigma)+1}.\]

A permutation statistic is called a Stirling statistic if
\begin{align}
&\sum_{\sigma \in \mathfrak{S}_{n}} x^{{\rm sst}(\sigma)}=x(x+1)(x+2)\cdots (x+n-1).
\end{align}
That is, $sst$ has the same generating function as the unsigned Stirling  number of the first kind. See \cite[Proposition 1.3.7]{Stanley-2012}. Here we describe   five Stirling statistics. The first is  the number of cycles in a decomposition of $\sigma$ into disjoint cycles, including those of length
1, which is denoted ${\rm cyc}(\sigma)$. For   the permutation $\sigma=2\,7\,1\,8\,3\,6\,5\,4$, its cycle decomposition is $(1\,2\,7\,5)\, (4\,8)\, (6)$, and so ${\rm cyc}(\sigma)=3$. Let $\sigma=\sigma_1\sigma_2\cdots \sigma_n \in \mathfrak{S}_n$. A left-to-right  maximum (resp. a left-to-right  minimum) of $\sigma$ is an element $\sigma_i$ such that $\sigma_{j}<\sigma_i$ (resp. $\sigma_{j}>\sigma_i$) for every $j<i$ and a right-to-left maximum (resp. a right-to-left minimum) of $\sigma$ is an element $\sigma_i$ such that $\sigma_{j}<\sigma_i$ (reps. $\sigma_{j}>\sigma_i$ ) for every $j>i$. Let ${\rm LRmax}(\sigma)$,  ${\rm LRmin}(\sigma)$, ${\rm RLmax}(\sigma)$ and  ${\rm RLmin}(\sigma)$ denote the  number of left-to-right maxima,    left-to-right minima,  right-to-left maxima  and right-to-left minima  of $\sigma$, respectively. For   the permutation $\sigma=2\,7\,1\,8\,3\,6\,5\,4$, we see that
\[{\rm LRmax}(\sigma)=3,\, {\rm LRmin}(\sigma)=2,\,{\rm RLmax}(\sigma)=4,\,{\rm RLmin}(\sigma)=3.\]

It is well known that
\begin{align} \label{def-Stirling}
&\sum_{\sigma \in \mathfrak{S}_{n}} x^{{\rm cyc}(\sigma)}=\sum_{\sigma \in \mathfrak{S}_{n}} x^{{\rm LRmax}(\sigma)}=\sum_{\sigma \in \mathfrak{S}_{n}} x^{{\rm RLmax}(\sigma)}=\sum_{\sigma \in \mathfrak{S}_{n}} x^{{\rm LRmin}(\sigma)}=\sum_{\sigma \in \mathfrak{S}_{n}} x^{{\rm RLmin}(\sigma)}\\[5pt]
&\quad \quad =x(x+1)(x+2)\cdots (x+n-1).
\end{align}

Carlitz and Scoville \cite{Carlitz-Scoville-1974}  considered the
following polynomials involving Eulerian-Stirling statistics, which we
refer to as the $(\alpha, \beta)$-Eulerian polynomials:
\begin{equation}\label{def-Euler-max}
A_n(x,y|{\alpha},{\beta})=\sum_{\sigma \in \mathfrak{S}_{n+1}} x^{{\rm asc}(\sigma)}y^{{\rm des}(\sigma)}{\alpha}^{{\rm LRmax}(\sigma)-1}{\beta}^{{\rm RLmax}(\sigma)-1}.
\end{equation}

They obtained the following generating function of $A_n(x,y|{\alpha},{\beta})$:
 \begin{thm}{\rm \!\! (Carlitz and Scoville \cite[Theorem 9]{Carlitz-Scoville-1974})}\label{C-S-main1}
 \begin{equation}\label{C-S-main1-eq}
 \sum_{n\geq 0} A_n(x,y|{\alpha},{\beta}) \frac{t^n}{n!}=\left(1+xF(x,y;t)\right)^{\alpha} \left(1+yF(x,y;t)\right)^{\beta},
 \end{equation}
 where $F(x,y;t)$ is given by
 \begin{equation}\label{defi-F}
 F(x,y;t)=\frac{e^{xt}-e^{yt}}{xe^{yt}-ye^{xt}}.
 \end{equation}
  \end{thm}
It's worth mentioning that  Carlitz and Scoville \cite{Carlitz-Scoville-1974}  used the terms fall for descent and rise for ascent. They referred to a left-to-right maximum and a right-to-left maximum of a permutation as a left upper record and a right upper record of a permutation.

Note that when $\alpha=0$,  ${\beta}=1$ and $x=1$, the polynomials $A_n(x,y|{\alpha},{\beta})$ reduce to the classical Eulerian polynomials $A_n(y)$ given by \eqref{def-Euler}. Accordingly,  we recover the generating function \eqref{gf-Euler} of $A_n(y)$ by setting $\alpha=0$,  ${\beta}=1$ and $x=1$ in \eqref{C-S-main1-eq}.

 When $x=y=1$ and $\beta=0$, it is not difficult to find that
 \begin{equation}\label{gf-stirling-special-t2}
 A_n(1,1|\alpha,0)=\sum_{{\sigma} \in \mathfrak{S}_{n}}{\alpha}^{{\rm LRmax}({\sigma})}
 \end{equation}
 and by \eqref{C-S-main1-eq}, we see that
  \begin{equation}\label{gf-stirling-special}
 \sum_{n\geq 0} A_n(1,1|\alpha,0) \frac{t^n}{n!}=\left(1+F(1,1;t)\right)^{\alpha}=\left(\frac{1}{1-t}\right)^\alpha.
 \end{equation}
Comparing the coefficients of $t^{n}/n!$ yields the generating function \eqref{def-Stirling}.

We would also like to mention that Foata and Sch\"utzenberger \cite{Foata-Schutzenberger-1970} introduced the following polynomial, which incorporates Eulerian-Stirling statistics  and has been referred to as the $q$-analogue of Eulerian polynomials by Brenti \cite{Brenti-2000}:
 \begin{equation}\label{Fota-Schutz}
  A_n(x,q)=\sum_{\sigma \in \mathfrak{S}_n}x^{{\rm exc}(\sigma)} q^{{\rm cyc}(\sigma)}.
  \end{equation}
  The grammar of $A_n(x,q)$ was found by Ma-Ma-Yeh-Zhu \cite{Ma-Ma-Yeh-Zhu-2018}.
 Brenti  \cite{Brenti-2000} obtained  the following generating function of $A_n(x,q)$ and showed that $A_n(x,q)$ is log-concave and unimodal.
 \begin{equation}\label{gf-brenti}
 1+\sum_{n\geq 1} A_{n}(x,q)\frac{t^n}{n!}=\left(\frac{e^{t(x-1)}-x}{1-x}\right)^{-q}.
\end{equation}
Thanks to the first fundamental transformation of Foata and Sch\"utzenberger \cite{Foata-Schutzenberger-1970}, we see that $q$-Eulerian polynomial \eqref{Fota-Schutz} is a special case of $(\alpha,\beta)$-Eulerian polynomials $A_n(x,y|{\alpha},{\beta})$. More precisely, we have
\[A_n(x,q)=\sum_{\sigma \in \mathfrak{S}_n}x^{{\rm exc}(\sigma)} q^{{\rm cyc}(\sigma)}=\sum_{\overline{\sigma} \in \mathfrak{S}_{n}}x^{{\rm des}(\overline{\sigma})}{q}^{{\rm LRmax}(\overline{\sigma})}=A_n(1,x|q,0).\]
Consequently,  we can retrieve the generating function \eqref{gf-brenti} for $A_n(x,q)$ by setting $x=1$ and $\beta=0$ and replacing $y$ and $\alpha$ with $x$ and $q$ respectively in \eqref{C-S-main1-eq}. It should be noted that the polynomials $A_n(x,q)$ has been shown to be related to the $1/k$-Eulerian polynomials, see Savage and Viswanathan \cite{Savage-Viswanathan-2012} and Ma and Mansour \cite{Ma-Mansour-2015} for example.

 The descent statistics are related to the Eulerian statistics, which are permutation statistics that depend only on the descent and length of a permutation, see  Gessel and Zhuang\cite{Gessel-Zhuang-2018}, Zhuang \cite{Zhuang-2016, Zhuang-2017}.   The classical descent statistics include variations of peaks and valleys, double ascents and double descents. In this paper, we adopt the
 terminology for these descent statistics provided by Gessel and Zhuang \cite{Gessel-Zhuang-2018} and  Zhuang \cite{Zhuang-2016, Zhuang-2017}.
 For the detailed definitions of these descent statistics,  see Section 2.

 Carlitz and Scoville \cite{Carlitz-Scoville-1974}   considered the generating function on the joint distribution of the number of exterior peaks, the number of  descents and the number of  ascents.
It should be noted that an exterior peak is referred to as a maximum by Carlitz and Scoville \cite{Carlitz-Scoville-1974}.

 \begin{thm}{\rm \!\! (Carlitz and Scoville \cite[Theorem 2]{Carlitz-Scoville-1974})}  \label{C-S-main2}
 Let $W(\sigma)$ denote  the number of exterior peaks of $\sigma$ (see Definition \ref{Dvp}).  Then
  \begin{equation}\label{C-S-main2-eq}
 \sum_{n\geq 0}\left(\sum_{\sigma \in \mathfrak{S}_{n+1}}u^{{\rm W}(\sigma)-1}v^{{\rm des}(\sigma)}w^{{\rm asc}(\sigma)}\right)\frac{t^n}{n!}=\left(1+yF(x,y;t)\right)
 \left(1+xF(x,y;t)\right),
 \end{equation}
 where $F(x,y;t)$ was given by  \eqref{defi-F} and
\begin{equation}\label{defi-x}
x=\frac{(w+v)+\sqrt{(w+v)^2-4uvw}}{2}, \quad
y=\frac{(w+v)-\sqrt{(w+v)^2-4uvw}}{2}.
\end{equation}
\end{thm}
Goulden and Jackson \cite[Exercise 3.3.46]{Goulden-Jackson-2004} and Stanley \cite[Exercise 1.61]{Stanley-2012}  reformulated  Carlitz and Scoville's result in terms of left double ascents and right double descents, see Definition \ref{Dvda} and Definition \ref{Dvdd}.   It should be noted that Goulden and Jackson \cite[Exercise 3.3.46]{Goulden-Jackson-2004} and Stanley \cite[Exercise 1.61]{Stanley-2012} refer to a left double ascent as a double rise, a right double descent as a double fall. A valley is called a modified minimum by  Goulden and Jackson \cite[Exercise 3.3.46]{Goulden-Jackson-2004}.

Let $V(\sigma)$, ${\rm lda}(\sigma)$ and ${\rm rdd}(\sigma)$ denote  the number of valleys of $\sigma$, left double ascents of $\sigma$, right double descents of $\sigma$ respectively (see Section 2). Goulden and Jackson \cite[Exercise 3.3.46]{Goulden-Jackson-2004} and Stanley \cite[Exercise 1.61]{Stanley-2012}
reformulated  Theorem \ref{C-S-main2} as
\begin{equation}\label{C-S-main2-varia-orign}
 \sum_{n\geq 1}\left(\sum_{\sigma \in \mathfrak{S}_{n}}u_1^{{\rm V}(\sigma)} u_2^{{\rm W}(\sigma)-1}u_3^{{\rm lda}(\sigma)}u_4^{{\rm rdd}(\sigma)}\right)\frac{t^n}{n!}=F(x,y;t),
 \end{equation}
 where $x+y=u_3+u_4$ and $xy=u_1u_2$.

Fu \cite{Fu-2018} provided a grammatical proof of \eqref{C-S-main2-varia-orign}.  Pan and Zeng \cite{Pan-Zeng-2019} provided {\it inv} $q$-analogue of \eqref{C-S-main2-varia-orign}.  Differentiating  both sides of \eqref{C-S-main2-varia-orign} with respect to $t$, we have

\begin{thm}{\rm \!\! (Carlitz and Scoville II \cite[Theorem 2]{Carlitz-Scoville-1974})}  \label{C-S-main2-varia}
  \begin{equation} \label{C-S-main2-varia-eq}
 \sum_{n\geq 0}\left(\sum_{\sigma \in \mathfrak{S}_{n+1}}u_1^{{\rm V}(\sigma)} u_2^{{\rm W}(\sigma)-1}u_3^{{\rm lda}(\sigma)}u_4^{{\rm rdd}(\sigma)}\right)\frac{t^n}{n!}=\left(1+yF(x,y;t)\right)
 \left(1+xF(x,y;t)\right),
 \end{equation}
 where $x+y=u_3+u_4$ and $xy=u_1u_2$.
\end{thm}
Setting $u_1=u_3=w$, $u_2=uv$ and $u_4=v$ in Theorem \ref{C-S-main2-varia}, and using  \eqref{rel-despkdd-b} and \eqref{rel-despkdd-a}, one could recover Theorem \ref{C-S-main2}.

The main objective of this paper is to investigate the following polynomial  involving descent-Stirling statistics.
\begin{equation}\label{defi-P}
    P_n(u_1,u_2,u_3,u_4|{\alpha},{\beta})=\sum_{\sigma \in \mathfrak{S}_{n+1}}u_1^{{\rm V}(\sigma)}u_2^{{\rm W}(\sigma)-1}u_3^{{\rm rdd}(\sigma)}u_4^{{\rm lda}(\sigma)}{\alpha}^{{\rm LRmax}(\sigma)-1}{\beta}^{{\rm RLmax}(\sigma)-1}.
\end{equation}
By using  the grammatical calculus introduced by Chen \cite{Chen-1993}, we establish the connection between the generating function of $P_n(u_1,u_2,u_3,u_4|{\alpha},{\beta})$ and the generating function \eqref{def-Euler-max} of $A_n(x,y|{\alpha},{\beta})$.

\begin{thm} We have\label{main2}
 \begin{equation}
 \sum_{n\geq 0}P_n(u_1,u_2,u_3,u_4|{\alpha},{\beta})\frac{t^n}{n!}
 =\left(1+yF(x,y;t)\right)^\frac{{\alpha}+{\beta}}{2}\left(1+xF(x,y;t)\right)^\frac{{\alpha}+{\beta}}{2}e^{\frac{1}{2}({\beta}-{\alpha})(u_3-u_4)t},
 \end{equation}
 where $x+y=u_3+u_4$, $xy=u_1u_2$ and $F(x,y;t)$ is given by \eqref{defi-F}.
\end{thm}

When ${\alpha}={\beta}=1$ in Theorem \ref{main2},   we  recover Theorem \ref{C-S-main2-varia}.

Using Theorem \ref{main2}, we derive the following   generating function of $P_n(u_1,u_2,u_3,u_4|{\alpha},{\beta})$.

\begin{thm} \label{main2c} We have
\begin{align*}
&\sum_{n\geq 0} P_n(u_1,u_2,u_3,u_4|{\alpha},{\beta}) \frac{t^n}{n!}=e^{\frac{1}{2}({\beta}-{\alpha})(u_3-u_4)t}\times \left(\cosh\left(\frac{t}{2}\sqrt{(u_3+u_4)^2  -4u_1u_2}\right) \right.\\[5pt]
& \hskip 2.5cm  \left. -\frac{u_3+u_4}{\sqrt{(u_3+u_4)^2-4u_1u_2}}\sinh\left(\frac{t}{2}\sqrt{(u_3+u_4)^2-4u_1u_2}\right)\right)^{-({\alpha}+{\beta})}.
\end{align*}
\end{thm}
As  applications of Theorem \ref{main2c}, we obtain $(\alpha,\beta)$-extensions of the   generating functions of peaks, left peaks, double  ascents, right double  ascents and left-right double ascents (see Theorems  \ref{gf-pk-thm},  \ref{gf-lrda}, \ref{gf-da},  \ref{gf-pk-lr} and \ref{gf-rda}). For more detailed explanations  on the   generating functions of these statistics, please see Section 2.

Based on Theorem \ref{main2c}, we obtain the following explicit expression of $P_n(u_1,u_2,u_3,u_4;{\alpha},{\beta})$ when ${\alpha}+{\beta}=-1$. In particular, we obtain the following interesting enumerative consequences.

\begin{thm}\label{enumer-pk} Let $M(\sigma)$ denote the number of interior peaks of $\sigma$.  When ${\alpha}+{\beta}=-1$ and for $n\geq 1$,
\begin{align}
&\sum_{\sigma \in \mathfrak{S}_{n+1}}u^{  {\rm M}(\sigma)} {\alpha}^{{\rm LRmin}(\sigma)-1}{\beta}^{{\rm RLmin}(\sigma)-1}  =\begin{cases} (1-u)^{\lfloor \frac{n}{2}\rfloor},  \quad \text{if \ $n$ is even},  \\[5pt]
-(1-u)^{\lfloor \frac{n}{2}\rfloor},  \quad \text{if \ $n$ is odd}.
\end{cases}
\end{align}
\end{thm}
\begin{thm} \label{enumer-lpk} Let  $L(\sigma)$ denote the number of left peaks of $\sigma$. For $n\geq 1$,
\[\sum_{\sigma \in \mathfrak{S}_{n}}u^{  {\rm L}(\sigma)}  (-1)^{{\rm RLmin}(\sigma)}=\begin{cases} (1-u)^{\lfloor \frac{n}{2}\rfloor},  \quad \text{if \ $n$ is even},  \\[5pt]
-(1-u)^{\lfloor \frac{n}{2}\rfloor},  \quad \text{if \ $n$ is odd}.
\end{cases}
\]
\end{thm}
\begin{thm}\label{enumer-rda} Let  $ {\rm rda}(\sigma)$ denote the number of right double ascents of $\sigma$. For $n\geq 1$,
\begin{align}
&2^n\sum_{\sigma \in \mathfrak{S}_{n+1}}u^{  {\rm rda}(\sigma)}  \left(-\frac{1}{2}\right)^{{\rm LRmin}(\sigma)+{\rm RLmin}(\sigma)-2} \nonumber \\[5pt]
&\quad \quad  =\begin{cases} ((1+u)^2-4)^{\lfloor \frac{n}{2}\rfloor},  \quad \text{if \ $n$ is even},  \\[5pt]
-(1+u)((1+u)^2-4)^{\lfloor \frac{n}{2}\rfloor},  \quad \text{if \ $n$ is odd}.
\end{cases}
\end{align}
\end{thm}
 Combining  Theorem \ref{C-S-main1} and Theorem \ref{main2},   we derive  two relations between $P_n(u_1,u_2,\break$
 $u_3,u_4|{\alpha},{\beta})$ and $A_n\left(x,y|\alpha, \beta\right)$ (see Theorem \ref{main2-cora} and Theorem \ref{main2-corb}).
These two relations enable us to establish  $(\alpha,\beta)$-extensions of the relations related to the Eulerian polynomial due to Stembridge, Petersen, Br\"and\'en and Zhuang, see Theorems \ref{main2-cora-1}, \ref{zhuang-pk-ex},  \ref{rela-lpk-exte} and \ref{rela-zhuang-lpk-exte}. We also obtain the  left peak version of Stembridge's formula and peak version of Petersen's formula (see Theorem \ref{leftpk-Stembridge} and Theorem  \ref{pk-Peterson}) and  their $(\alpha,\beta)$-extensions, see Theorem \ref{main2-cora-2a} and Theorem  \ref{main-coraa}. The following two consequences can be viewed as the $(\alpha,\beta)$-extensions of Stembridge's formula and the left peak version of Stembridge's formula
\begin{thm}\label{main2-cora-1}For $n\geq 1$,
\begin{align}
&\sum_{\sigma \in \mathfrak{S}_{n}}(xy)^{{\rm M}(\sigma)}\left(\frac{x+y}{2}\right)^{n-2{\rm M}(\sigma)-1}{\alpha}^{{\rm LRmin}(\sigma)-1}{\beta}^{{\rm RLmin}(\sigma)-1}\nonumber \\[5pt]
&=\sum_{\sigma \in \mathfrak{S}_{n}}x^{{\rm des}(\sigma)}y^{n-{\rm des}(\sigma)-1}{\left(\frac{{\alpha}+{\beta}}{2}\right)}^{{\rm LRmin}(\sigma)+{\rm RLmin}(\sigma)-2},
\end{align}
where $M(\sigma)$ counts the number of interior peaks of $\sigma$.
\end{thm}

\begin{thm} \label{main2-cora-2a} For $n\geq 0$,
\begin{align}\label{C-S-3-varia-g}
&\sum_{\sigma \in \mathfrak{S}_{n}}(xy)^{{\rm L}(\sigma)}\left(\frac{x+y}{2}\right)^{n-2{\rm L}(\sigma)}{\beta}^{{\rm RLmin}(\sigma)}\nonumber \\[5pt]
&=\sum_{\sigma \in \mathfrak{S}_{n+1}}x^{{\rm des}(\sigma)}y^{n-{\rm des}(\sigma)} {\left(\frac{{\beta}}{2}\right)}^{{\rm LRmin}(\sigma)+{\rm RLmin}(\sigma)-2},
\end{align}
where $L(\sigma)$ counts the number of left peaks of $\sigma$.
\end{thm}
Specializing Theorem \ref{main2-cora-1} and Theorem \ref{main2-cora-2a} allows us to derive the $(\alpha,\beta)$-extensions of the tangent and secant numbers.  Recall that the tangent  number $E_{2n+1}$ and the secant  number $E_{2n}$ are defined by
\[\sum_{n\geq 0} E_{2n+1}\frac{t^{2n+1}}{(2n+1)!}=\tan(t) \quad \text{and}  \quad \sum_{n\geq 0} E_{2n}\frac{t^{2n}}{(2n)!}=\sec(t). \]
A permutation $\sigma=\sigma_1\sigma_2\cdots \sigma_n \in \mathfrak{S}_n$  is down-up  (or alternating) if $\sigma_1 > \sigma_2 < \sigma_3 > \sigma_4 <\ldots$ and a permutation $\sigma=\sigma_1\sigma_2\cdots \sigma_n \in \mathfrak{S}_n$  is up-down  (or reverse alternating) if $\sigma_1 < \sigma_2 > \sigma_3 < \sigma_4 >\ldots$.   The down-up permutations in $\mathfrak{S}_4$ are $2\,1\,4\,3$, $3\,1\,4\,2$, $3\,2\,4\,1$, $4\,1\,3\,2$, $4\,2\,3\,1$ and the up-down permutations in $\mathfrak{S}_4$ are $3\,4\,1\,2$, $2\,4\,1\,3$, $2\,3\,1\,4$, $1\,4\,2\,3$, $1\,3\,2\,4$. It is easy to show that the number of  down-up permutations of $[n]$ equals the number of up-down permutations of $[n]$. Andr\'e \cite{Andre-1879} showed that $E_n$ counts the number of down-up (or up-down) permutations of $[n]$.

Euler \cite{Euler-1755}   found the following interesting relation: For $n\geq 1$,
\begin{align}\label{relation-Euler}
&\sum_{\sigma \in \mathfrak{S}_{n}}(-1)^{{\rm exc}(\sigma)}=\begin{cases}
 (-1)^{\frac{n-1}{2}} E_n, &\text{ if } n \text{ is odd,} \\[5pt]
0, &\text{ if } n \text{ is even}. \\[5pt]
 \end{cases}
 \end{align}
Roselle \cite{Roselle-1968}   obtained the following parallel result to Euler involving secant numbers:  For $n\geq 1$,
  \begin{align}\label{relation-Euler-Roselle}
&\sum_{\sigma \in \mathfrak{D}_{n}}(-1)^{{\rm exc}(\sigma)}=\begin{cases}
 (-1)^{\frac{n}{2}} E_n, &\text{ if } n \text{ is even,} \\[5pt]
0, &\text{ if } n \text{ is odd}, \\[5pt]
 \end{cases}
 \end{align}
  where $\mathfrak{D}_{n}$ counts the number of permutations in $\mathfrak{S}_n$ without fixed points. It should be noted many different $q$-analogues of \eqref{relation-Euler} and \eqref{relation-Euler-Roselle} have been established by \cite{Foata-Han-2010, Josuat-Verges-2010, Shin-Zeng-2010a, Shin-Zeng-2010b}.

 Setting $x=-1$ and $y=1$ in Theorem  \ref{main2-cora-1},   we have
 \begin{thm}\label{main2-coraa-alternate}Let $\mathfrak{S}^a_n$ denote the set of up-down permutations of $[n]$.   For $n\geq 1$,
 \begin{align}
&\sum_{\sigma \in \mathfrak{S}_{n}}(-1)^{{\rm des}(\sigma)}{\left(\frac{{\alpha}+{\beta}}{2}\right)}^{{\rm LRmin}(\sigma)+{\rm RLmin}(\sigma)-2}\nonumber\\[5pt]
&=\begin{cases}
 (-1)^{\frac{n-1}{2}} \sum_{\sigma \in \mathfrak{S}^a_{n}} {\alpha}^{{\rm LRmin}(\sigma)-1}{\beta}^{{\rm RLmin}(\sigma)-1}, &\text{ if } n \text{ is odd,} \\[5pt]
0, &\text{ if } n \text{ is even}. \\[5pt]
 \end{cases}
 \end{align}
 \end{thm}

 Setting $x=-1$ and $y=1$ in   Theorem \ref{main2-cora-2a}, we have
\begin{thm}\label{main2-coraa-alternate2}Let $\mathfrak{S}^a_n$ denote the set of down-up permutations of $[n]$.   For $n\geq 1$,
 \begin{align}
&\sum_{\sigma \in \mathfrak{S}_{n+1}}(-1)^{{\rm des}(\sigma)}{\left(\frac{{\beta}}{2}\right)}^{{\rm LRmin}(\sigma)+{\rm RLmin}(\sigma)-2}\nonumber\\[5pt]
&=\begin{cases}
 (-1)^{\frac{n}{2}} \sum_{\sigma \in \mathfrak{S}^a_{n}} {\beta}^{{\rm RLmin}(\sigma)}, &\text{ if } n \text{ is even,} \\[5pt]
0, &\text{ if } n \text{ is odd}. \\[5pt]
 \end{cases}
 \end{align}
 \end{thm}

Setting $\alpha=\beta=1$ in Theorem \ref{main2-coraa-alternate}, we could recover   Euler's relation \eqref{relation-Euler} with the aid of the first fundamental transformation of Foata and Sch\"utzenberger \cite{Foata-Schutzenberger-1970}.
Setting $\beta=1$ in Theorem \ref{main2-coraa-alternate2}, we obtain the following identity, which seems to be new.
 \begin{align}\label{relation-Euler-Roselle-varia}
&\sum_{\sigma \in \mathfrak{S}_{n+1}}(-1)^{{\rm des}(\sigma)}{\left(\frac{1}{2}\right)}^{{\rm LRmin}(\sigma)+{\rm RLmin}(\sigma)-2}=\begin{cases}
 (-1)^{\frac{n}{2}} E_n, &\text{ if } n \text{ is even,} \\[5pt]
0, &\text{ if } n \text{ is odd}. \\[5pt]
 \end{cases}
 \end{align}

 Combining \eqref{relation-Euler-Roselle} and \eqref{relation-Euler-Roselle-varia}, we obtain the following identity:
\begin{equation}
\sum_{\sigma \in \mathfrak{S}_{n+1}}(-1)^{{\rm des}(\sigma)}{\left(\frac{1}{2}\right)}^{{\rm LRmin}(\sigma)+{\rm RLmin}(\sigma)-2}=\sum_{\sigma \in \mathfrak{D}_{n}}(-1)^{{\rm exc}(\sigma)}.
\end{equation}
 It would be interesting to give a combinatorial proof of the above identity.

This paper is organized as follows.  Section 2 provides a review of   some classical descent statistics, including left peaks, interior peaks, exterior peaks, valleys, left double ascents, double ascents, right double ascents, left-right ascents, left double descents, double descents, right double descents, left-right descents. We then collect the generating functions associated with these statistics and their relations with the Eulerian polynomials. Section 3 is dedicated to proving the main result of this paper (Theorem \ref{main2}) using the grammatical calculus introduced by Chen \cite{Chen-1993}.  In Section 4, we first   derive the generating function of $P_n(u_1,u_2,u_3,u_4|{\alpha},{\beta})$   using Theorem \ref{main2}. We then present  $(\alpha,\beta)$-extensions of  some  known generating functions related to descent statistics by specializing the generating function of $P_n(u_1,u_2,u_3,u_4|{\alpha},{\beta})$.  Section 5 aims to establish an explicit expression of $P_n(u_1,u_2,u_3,u_4|{\alpha},{\beta})$ when $\alpha+\beta=-1$. This result can be specialized to obtain Theorems \ref{enumer-pk}, \ref{enumer-lpk}, and \ref{enumer-rda}. In Section 6,
we first establish two relations between $P_n(u_1,u_2,u_3,u_4|{\alpha},{\beta})$ and $A_n\left(x,y|\alpha, \beta\right)$  using Theorem \ref{main2}. We then derive $(\alpha,\beta)$-extensions of some known relations between descent statistics and the Eulerian polynomials given by Stembridge, Petersen, Br\"and\'en and Zhuang  by specializing these two relations.

\section{Descent statistics}

In this section, we begin by revisiting classical descent statistics, which encompass variations related to peaks, valleys, double ascents, and double descents. We then collect the generating functions associated with these statistics and   their relationships with Eulerian polynomials. Here we follow Stanley's terminology for peaks and its variations, as described in \cite[Exercise 1.61]{Stanley-2012}. For double ascents, double descents and their variations, we adhere to the definitions provided by Zhuang in \cite{Zhuang-2016}.

 \begin{defi}[Variations in peaks]\label{Dvp}  Given a permutation $\sigma=\sigma_1\cdots \sigma_n$,
 \begin{itemize}
\item[{\rm (1)}] we say that  $i$ is a {\bf left peak} of $\sigma$ if $1\leq i<n$ and $\sigma_{i-1}<\sigma_i>\sigma_{i+1}$ under the assumption that $\sigma_0=\sigma_{n+1}=0$. Let $L(\sigma)$ denote the number of left peaks of $\sigma$.

\item[{\rm (2)}] we say that   $i$ is   an {\bf
 interior   peak} (or a peak for short) of $\sigma$ if $1< i<n$ and $\sigma_{i-1}<\sigma_i>\sigma_{i+1}$.
Let $M(\sigma)$ denote the number of   peaks of $\sigma$.

\item[{\rm (3)}] we say that   $i$ is an {\bf  exterior peak} of $\sigma$ if $1\leq i\leq n$ and $\sigma_{i-1}<\sigma_i>\sigma_{i+1}$ under the assumption that $\sigma_0=\sigma_{n+1}=0$.
Let $W(\sigma)$ denote the number of exterior peaks of $\sigma$.

\item[{\rm (4)}] we say that $i$ is   a {\bf  valley}  of $\sigma$ if $1\leq i\leq n$  and  $\sigma_{i-1}>\sigma_{i}<\sigma_{i+1}$ under the assumption that $\sigma_0=\sigma_{n+1}=0$. Let   ${\rm V}(\sigma)$ denote  the number of valleys of $\sigma$.
\end{itemize}

 \end{defi}

 For  the permutation $\sigma=7\,1\,3\,8\,5\,9\,6\,2\,4 \in \mathfrak{S}_9$,
 we see that
 \[L(\sigma)=3,\quad  M(\sigma)=2, \quad W(\sigma)=4,\quad {\rm V}(\sigma)=3. \]
Note that the symbols $L(\sigma)$, $M(\sigma)$, and $W(\sigma)$ used in this context were introduced by Chen and Fu \cite{Chen-Fu-2022},  which are   meaningful and easy to remember. The letter $L$ looks like having a peak on the left. It should also be noted that Carlitz and Scoville \cite{Carlitz-Scoville-1974} refer to an exterior peak as a maximum, while Goulden and Jackson \cite[Exercise 3.3.46]{Goulden-Jackson-2004} describe it as a modified maximum. Similarly, Goulden and Jackson \cite[Exercise 3.3.46]{Goulden-Jackson-2004} term a valley as a modified minimum.

It is evident that for $\sigma \in \mathfrak{S}_n$,
  \begin{equation}\label{rel-despkdd-b}
 {\rm V}(\sigma)={\rm W}(\sigma)-1
\end{equation}
and
\begin{equation}\label{rel-despkdd-b}
 \sum_{{\sigma \in \mathfrak{S}_{n}}}u^{  {\rm M}(\sigma)}= \sum_{{\sigma \in \mathfrak{S}_{n}}}u^{  {\rm V}(\sigma)}=\sum_{{\sigma \in \mathfrak{S}_{n}}}u^{ {\rm W}(\sigma)-1}
\end{equation}
 Hence it suffices to consider left peak and interior peak.  The  generating function for the left peak polynomials is attribute to  Gessel \cite[Sequence A008971]{OLEN}, see Zhuang \cite[Theorem 10]{Zhuang-2016}.
\begin{thm}[Gessel] We have
\begin{align}\label{gf-lpk}
\sum_{n\geq 0}\left( \sum_{\sigma \in \mathfrak{S}_{n}}u^{  {\rm L}(\sigma)}\right) \frac{t^n}{n!} =  \frac{\sqrt{1-u}}{\sqrt{1-u} \cosh\left(t\sqrt{1-u}\right)-\sinh\left(t\sqrt{1-u}\right)}.
\end{align}
\end{thm}
The following presents  the generating function for the peak polynomials.  As brought up by Stanley \cite{Stanley-2012}, the generating function of the peak polynomials can be deduced from an equation of David-Barton \cite{David-Barton-1962}, see Chen and Fu \cite{Chen-Fu-2022-a} for more information.  The equivalent formulae have been found by Entringer \cite{Entringer-1969},  Kitaev \cite{Kitaev-2007} and Zhuang  \cite[Theorem 9]{Zhuang-2016}.
\begin{thm} We have
\begin{align}\label{gf-pk}
&\sum_{n\geq 0}\left( \sum_{{\sigma \in \mathfrak{S}_{n}}}u^{  {\rm M}(\sigma)}\right) \frac{t^n}{n!} =  \frac{\sqrt{1-u}\cosh\left(t\sqrt{1-u}\right)}{\sqrt{1-u} \cosh\left(t\sqrt{1-u}\right)-v\sinh\left(t\sqrt{1-u}\right)}.
\end{align}
\end{thm}

Chen and Fu \cite{Chen-Fu-2023} provide a grammatical proof of \eqref{gf-lpk} and \eqref{gf-pk}.

 Stembridge \cite{Stembridge-1997} first considered the relation between the peak polynomials and the Eulerian polynomials. He  obtained the following relation between the peak polynomials and  the Eulerian polynomials in the study of his theory of enriched $P$-partitions, which was rediscovered by Br\"and\'en \cite{Branden-2008}  with the aid of  the ``modified Foata-Strehl action", a variant of a group action on permutations originally defined by Foata and Strehl \cite{Foata-Strehl-1974}.

\begin{thm}[Stembridge] For $n\geq 1$,
\begin{equation}\label{Stembridge-v}
 \sum_{\sigma \in \mathfrak{S}_{n}} x^{{\rm des}(\sigma)}=\left(\frac{1+x}{2}\right)^{n-1} \sum_{\sigma \in \mathfrak{S}_{n}}\left(\frac{4x}{(1+x)^2}\right)^{{\rm M}(\sigma)}.
 \end{equation}
 \end{thm}

  Petersen \cite[Observation 3.1.2]{Petersen-2006} established a relation between the left peak polynomials  and the Eulerian  polynomials, stated as follows.
\begin{thm}[Petersen]
For $n\geq 1$,
\begin{equation}\label{relation-lpk}
(1+x)^n\sum_{\sigma \in \mathfrak{S}_{n}}\left(\frac{4x}{(1+x)^2}\right)^{{\rm L}(\sigma)}=\sum_{k=1}^n {n\choose k} 2^k(1-x)^{n-k}\sum_{\sigma \in \mathfrak{S}_{k}}x^{{\rm des}(\sigma)+1}+(1-x)^n.
\end{equation}
\end{thm}

Chen and Fu \cite{Chen-Fu-2022-a}  provided   grammatical proofs of \eqref{Stembridge-v} and \eqref{relation-lpk}.

Recently, Zhuang \cite{Zhuang-2017} established two relations between the joint polynomials of peaks (or left peaks ) and descents  and the Eulerian  polynomials.

\begin{thm} {\rm (Zhuang \cite[Theorem 4.2]{Zhuang-2017})}
For $n\geq 1$,
\begin{equation}\label{zhuang-pk-des}
\sum_{\sigma \in \mathfrak{S}_n} u^{{\rm M}(\sigma)+1}v^{{\rm des}(\sigma)+1}=\left(\frac{1+b}{1+ab}\right)^{n+1}\sum_{\sigma \in \mathfrak{S}_n} a^{{\rm des}(\sigma)+1},
\end{equation}
where
\begin{equation} \label{zhuang-relb}
a=\frac{(1+v)^2-2uv-(1+v)\sqrt{(1+v)^2-4uv}}{2uv}
\end{equation}
and
\begin{equation}\label{zhuang-rela}
b=\frac{1+v^2-2uv-(1-v)\sqrt{(1+v)^2-4uv}}{2(1-u)v}.
\end{equation}
\end{thm}

\begin{thm}{\rm (Zhuang \cite[Theorem 4.7]{Zhuang-2017})}
For $n\geq 1$,
\begin{equation}\label{zhuang-lpk-des}
\sum_{\sigma \in \mathfrak{S}_n} u^{{\rm L}(\sigma)}v^{{\rm des}(\sigma)}=\frac{1}{(1+ab)^n}\left(\sum^n_{k=1}{n\choose k}(1+b)^k(1-a)^{n-k}\sum_{\sigma \in \mathfrak{S}_k} a^{{\rm des}(\sigma)+1}+(1-a)^n\right),
\end{equation}
where $a$ and $b$ are defined by \eqref{zhuang-relb} and \eqref{zhuang-rela}.
\end{thm}

 \begin{defi}[Variations in double ascents] \label{Dvda} Given a permutation $\sigma=\sigma_1\cdots \sigma_n$,
 \begin{itemize}
\item[{\rm (1)}] we say that  $i$ is  a {\bf left double ascent} of $\sigma$ if $1\leq i<n$  and  $\sigma_{i-1}<\sigma_{i}<\sigma_{i+1}$ under the assumption that
$\sigma_0=0$.  Let ${\rm {lda}}(\sigma)$ denote the  number of left double ascents of $\sigma$.

\item[{\rm (2)}] we say that  $i$ is   a {\bf double ascent} of $\sigma$   if $1< i<n$  and  $\sigma_{i-1}<\sigma_{i}<\sigma_{i+1}$.  Let ${\rm {da}}(\sigma)$ denote the  number of  double ascents of $\sigma$.

\item[{\rm (3)}] we say that $i$ is    a {\bf right double ascent} of $\sigma$ if $1< i\leq n$  and  $\sigma_{i-1}<\sigma_{i}<\sigma_{i+1}$ under the assumption that
$\sigma_{n+1}=+\infty$.  Let ${\rm {rda}}(\sigma)$ denote the  number of right double ascents of $\sigma$.

\item[{\rm (4)}] we say that $i$ is     a {\bf left-right  double ascent} of $\sigma$ if   $1\leq i\leq n$  and  $\sigma_{i-1}<\sigma_{i}<\sigma_{i+1}$ under the assumption that
$\sigma_0=0$ and $\sigma_{n+1}=+\infty$. Let ${\rm {lrda}}(\sigma)$ denote the  number of  left-right double ascents of $\sigma$.
\end{itemize}

 \end{defi}

 \begin{defi}[Variations in double descents] \label{Dvdd} Given a permutation $\sigma=\sigma_1\cdots \sigma_n$,
 \begin{itemize}
\item[{\rm (1)}] we say that  $i$ is  a {\bf left double descent} of $\sigma$ if $1\leq i<n$  and  $\sigma_{i-1}>\sigma_{i}>\sigma_{i+1}$ under the assumption that
$\sigma_0=+\infty$.  Let ${\rm {ldd}}(\sigma)$ denote the  number of left double descents of $\sigma$.

\item[{\rm (2)}] we say that  $i$ is   a {\bf double descent} of $\sigma$   if $1< i<n$  and  $\sigma_{i-1}>\sigma_{i}>\sigma_{i+1}$.  Let ${\rm {dd}}(\sigma)$ denote the  number of  double descents of $\sigma$.

\item[{\rm (3)}] we say that $i$ is    a {\bf right double descent} of $\sigma$ if $1< i\leq n$  and  $\sigma_{i-1}>\sigma_{i}>\sigma_{i+1}$ under the assumption that
$\sigma_{n+1}=0$.  Let ${\rm {rdd}}(\sigma)$ denote the  number of right double descents of $\sigma$.

\item[{\rm (4)}] we say that $i$ is     a {\bf left-right  double ascent} of $\sigma$ if   $1\leq i\leq n$  and  $\sigma_{i-1}>\sigma_{i}>\sigma_{i+1}$ under the assumption that
$\sigma_0=+\infty$ and $\sigma_{n+1}=0$. Let ${\rm {lrdd}}(\sigma)$ denote the  number of  left-right double ascents of $\sigma$.
\end{itemize}

 \end{defi}

 For  the permutation $\sigma=7\,1\,3\,8\,5\,9\,6\,2\,4 \in \mathfrak{S}_9,$
 we see that
 \[{\rm {da}}(\sigma)=1,\quad  {\rm {lda}}(\sigma)=1, \quad {\rm {rda}}(\sigma)=2,\quad {\rm {lrda}}(\sigma)=2 \]
 and
  \[{\rm {dd}}(\sigma)=1,\quad  {\rm {ldd}}(\sigma)=2, \quad {\rm {rdd}}(\sigma)=1,\quad {\rm {lrdd}}(\sigma)=2. \]
By definition, we see that for $\sigma \in \mathfrak{S}_n$,
\begin{equation}\label{rel-despkdd-a}
 {\rm des}(\sigma)={\rm W}(\sigma)+{\rm rdd}(\sigma)-1,\quad \  \quad {\rm asc}(\sigma)={\rm W}(\sigma)+{\rm lda}(\sigma)-1
\end{equation}
and
\begin{equation}\label{rel-deslpkdd-a}
 {\rm des}(\sigma)={\rm L}(\sigma)+{\rm dd}(\sigma),\quad \  \quad {\rm asc}(\sigma)={\rm L}(\sigma)+{\rm lrda}(\sigma)-1.
\end{equation}
 It is evident from taking reverses and complements that we only need to consider double ascents, right double ascents and left-right double ascents.  By generalizing Gessel's reciprocity formula for noncommutative symmetric functions, Zhuang gave a systematic method for obtaining the generating functions for double ascents, right double ascents and left-right double ascents. It should be noted that the equivalent form of the generating function for double ascents was established by  Elizalde and Noy \cite{Elizalde-Noy-2003}.   Elizalde and Noy \cite{Elizalde-Noy-2003} referred to a double descent as a proper double descent.

\begin{thm}{\rm (\!\! Elizalde and Noy \cite{Elizalde-Noy-2003},  Zhuang \cite[Theorem 12]{Zhuang-2016})}
  \begin{align}\label{Zhuangda}
\sum_{n\geq 0}\left( \sum_{{\sigma \in \mathfrak{S}_{n}}}u^{ {\rm da}(\sigma)}  \right) \frac{t^n}{n!} =  \frac{v e^{\frac{(1-u)}{2}t}  }{{v \cosh\left(\frac{1}{2}vt \right)-(1+u)\sinh\left(\frac{1}{2}vt\right)}},
\end{align}
where $v=\sqrt{(u+1)^2-4}$.
\end{thm}

\begin{thm}{\rm (\!\! Zhuang \cite[Theorem 13]{Zhuang-2016})}
 \begin{align}\label{zhuang-rda}
\sum_{n\geq 0}\left( \sum_{{\sigma \in \mathfrak{S}_{n}}}u^{ {\rm rda}(\sigma)}  \right) \frac{t^n}{n!} =  \frac{v \cosh\left(\frac{1}{2}vt \right)+(1-u)\sinh\left(\frac{1}{2}vt\right)}{{v \cosh\left(\frac{1}{2}vt \right)-(1+u)\sinh\left(\frac{1}{2}vt\right)}}
\end{align}
and
 \begin{align}\label{zhuang-lrda}
\sum_{n\geq 0}\left( \sum_{{\sigma \in \mathfrak{S}_{n}}}u^{ {\rm lrda}(\sigma)}  \right) \frac{t^n}{n!} =  \frac{v e^{\frac{(u-1)}{2}t}  }{{v \cosh\left(\frac{1}{2}vt \right)-(1+u)\sinh\left(\frac{1}{2}vt\right)}},
\end{align}
where $v=\sqrt{(u+1)^2-4}$.
\end{thm}

\section{The grammatical  derivation for Theorem  \ref{main2}}

The main objective of this section is to  give a   proof of Theorem  \ref{main2} by using the grammatical calculus introduced by Chen \cite{Chen-1993}. A context-free grammar $G$ over a set $V=\{x,y,z,\ldots\}$ of variables is a set substitution rules replacing a variable in $V$ by a Laurent polynomial of variables in $V$. For a context-free grammar $G$ over $V$, the formal derivative $D$ with respect to $G$ is defined as a linear operator acting on Laurent polynomials with variables in $V$ such that each substitution rule is treated as the common differential rule that satisfies  the following relations:
\begin{align} \label{gramma-add-rela}
D(u+v)&=D(u)+D(v)\\[5pt]
D(uv)&=D(u)v+uD(v). \label{gramma-mult-rela}
\end{align}
Hence, it obeys the Leibniz's rule
\[D^{n}(uv)=\sum_{k=0}^n{n \choose k}D^k(u)D^{n-k}(v).
\]
For a constant $c$, we have $D(c)=0$.

A formal derivative $D$ with respect to $G$ is also associated with an exponential generating function. For a  Laurent polynomial $w$ of variables in $V$, let
\begin{equation}
    {\rm Gen}^{(G)}(w,t)=\sum_{n\geq 0}D^n(w)\frac{t^n}{n!}.
\end{equation}
Then, by \eqref{gramma-add-rela} and \eqref{gramma-mult-rela}, we derive that
\begin{equation}\label{gramma-add}
{\rm Gen}^{(G)}(u+v,t) = {\rm Gen}^{(G)}(u,t)+{\rm Gen}^{(G)}(v,t).
\end{equation}
\begin{equation}\label{gramma-multiple}
{\rm Gen}^{(G)}(uv,t) = {\rm Gen}^{(G)}(u,t){\rm Gen}^{(G)}(v,t).
\end{equation}
For more information on the grammatical calculus, we refer to Chen \cite{Chen-1993} and Chen and Fu \cite{Chen-Fu-2017, Chen-Fu-2022-a}.

Dumont \cite{Dumont-1996} showed the following grammar
\begin{equation}\label{C-S-G-orign}
G_1=\{x\rightarrow xy,\, y\rightarrow xy\}.
\end{equation}
generates the Eulerian polynomials $A_n(x)$. More precisely, let $D_{G_{1}}$ be the formal derivative with respect to the   grammar $G_1$ given by \eqref{C-S-G-orign}, then for $n\geq 1$,
\[D_{G_1}^n(y)=xA_{n}(x,y|0,1),\]
Here we adopt the notion $A_n(x,y|\alpha,\beta)$ given by \eqref{def-Euler-max} to represent  the bivariate  Eulerian polynomials, where
\[xA_{n}(x,y|0,1)=\sum_{\sigma \in \mathfrak{S}_{n}} x^{{\rm asc}(\sigma)+1}y^{{\rm des}(\sigma)+1}\]

  Chen and Fu  \cite{Chen-Fu-2022}  showed that
\begin{align}\label{chen-fu-eulera}
{\rm Gen}^{(G_1)}(y,t):=\sum_{n\geq 0}D_{G_1}^n(y)\frac{t^n}{n!} =x(1+yF(x,y;t)),
\end{align}
where $F(x,y;t)$ is given in \eqref{defi-F}. Together with Dumont's result, they provided a grammatical proof of the generating function \eqref{gf-Euler} of $A_n(x)$.

Similarly, it can be shown that
\begin{align}\label{chen-fu-eulerb}
{\rm Gen}^{(G_1)}(x,t):=\sum_{n\geq 0}D_{G_1}^n(x)\frac{t^n}{n!}&=
y(1+xF(x,y;t)).
\end{align}

In this section, we first show that the following grammar
  \begin{equation}\label{L-P-four}
    \widetilde{G}=\{a\rightarrow a\alpha u_4,\, b\rightarrow b\beta u_3,\,u_4\rightarrow u_1u_2, \, u_3\rightarrow u_1u_2,\,u_1\rightarrow u_1u_3, u_2\rightarrow  u_2u_4 \}.
     \end{equation}
 can be used to generate the polynomial $P_n(u_1,u_2,u_3,u_4|{\alpha},{\beta})$. More precisely,
\begin{thm}\label{thm-grammar-P}
Let  $D_{\widetilde{G}}$ be the formal derivative with respect to the   grammar defined in \eqref{L-P-four}, we have
\begin{equation}
    {D^{n}_{\widetilde{G}}}(ab)=abP_n
    (u_1,u_2,u_3,u_4|{\alpha},{\beta}).
\end{equation}
\end{thm}
Based on Theorem \ref{thm-grammar-P}, we give a proof of Theorem  \ref{main2} using the grammatical calculus. More precisely, it suffices to demonstrate the following theorem.
\begin{thm}\label{thm-grammar-P2}
Let  $D_{\widetilde{G}}$ be the formal derivative with respect to the   grammar defined in \eqref{L-P-four}, we have
\begin{equation}
    {\rm Gen}^{(\widetilde{G})}(ab,t)=ab\left(1+yF(x,y;t)\right)^\frac{{\alpha}+{\beta}}{2}
\left(1+xF(x,y;t)\right)^\frac{{\alpha}+{\beta}}{2}
e^{\frac{1}{2}({\beta}-{\alpha})(u_3-u_4)t},
\end{equation}
where $x+y=u_3+u_4$ and $xy=u_1u_2$ and $F(x,y;t)$ is given by \eqref{defi-F}.
\end{thm}

\subsection{A grammatical labeling of $P_n(u_1,u_2, u_3,u_4|{\alpha},{\beta})$}

To prove   Theorem \ref{thm-grammar-P}, we are required to give the combinatorial definition  of  $P_n(u_1,u_2, u_3,u_4|{\alpha},{\beta})$ involving the left-to-right minima and the right-to-left minima.  Recall that the complement of $\sigma=\sigma_1\cdots \sigma_n \in \mathfrak{S}_n$ is given by
\[\sigma^{c}=(n+1-\sigma_1)(n+1-\sigma_2)\cdots (n+1-\sigma_n).\]
For example, if $\sigma=7\,1\,3\,8\,5\,6\,2\,4$, then the complement of $\sigma$ is given by $\sigma^{c}=2\,8\,6\,1\,4\,3\,7\,5$.

Evidently, the complement provides a bijection between $\mathfrak{S}_n$ and $\mathfrak{S}_n$. Moreover, for $\sigma \in \mathfrak{S}_n$ and $\sigma^{c}$ is the complement of $\sigma$, we have
\begin{equation}\label{rel-compa}
{\rm asc}(\sigma)={\rm des}(\sigma^c),\quad {\rm des}(\sigma)={\rm asc}(\sigma^c),
\end{equation}
\begin{equation}\label{rel-compb}
{\rm lda}(\sigma)={\rm ldd}(\sigma^c),\quad {\rm rdd}(\sigma)={\rm rda}(\sigma^c),
\end{equation}
\begin{equation}\label{rel-compc}
{\rm LRmax}(\sigma)={\rm  LRmin}(\sigma^c),\quad {\rm RLmax}(\sigma)={\rm RLmin}(\sigma^c)
\end{equation}
and
\begin{equation}\label{rel-compd}
{\rm W}(\sigma)-1={\rm V}(\sigma)={\rm M}(\sigma^c).
\end{equation}
Hence we find that $A_n(x,y|{\alpha},
{\beta})$ and $P_n(u_1,u_2,u_3,u_4|{\alpha},{\beta})$ can also be interpreted as follows:

\begin{equation}\label{defi-A-min}
A_n(x,y|{\alpha},{\beta})=\sum_{\sigma \in \mathfrak{S}_{n+1}} x^{{\rm des}(\sigma)}y^{{\rm asc}(\sigma)}{\alpha}^{{\rm LRmin}(\sigma)-1}{\beta}^{{\rm RLmin}(\sigma)-1}
\end{equation}
and
\begin{equation}\label{defi-P-min}
       P_n(u_1,u_2,u_3,u_4|{\alpha},{\beta})=\sum_{\sigma \in \mathfrak{S}_{n+1}}(u_1u_2)^{{\rm M}(\sigma)} u_3^{{\rm rda}(\sigma)}u_4^{{\rm ldd}(\sigma)}{\alpha}^{{\rm LRmin}(\sigma)-1}{\beta}^{{\rm RLmin}(\sigma)-1}.
\end{equation}
Please take note that the polynomial $P_n(u_1,u_2,u_3,u_4|{\alpha},{\beta})$ encompasses interior peaks, right double ascents, and left double descents. However, it's worth mentioning that left peaks, left-right double ascents, and double descents can also be characterized by specializing the   polynomial $P_n(u_1,u_2,u_3,u_4|{\alpha},{\beta})$. More precisely, setting $\alpha=0$ in \eqref{defi-P-min}, we find that
\begin{align}\label{rda-lrdamin}
P_n(u_1,u_2,u_3,u_4|0,{\beta})&=\sum_{\sigma \in \mathfrak{S}_{n+1}}(u_1u_2)^{{\rm M}(\sigma)} u_3^{{\rm rda}(\sigma)}u_4^{{\rm ldd}(\sigma)}{0}^{{\rm LRmin}(\sigma)-1}{\beta}^{{\rm RLmin}(\sigma)-1}\nonumber \\[5pt]
&=\sum_{\substack{\sigma \in \mathfrak{S}_{n+1}\\[2pt]\sigma_{1}=1}}(u_1u_2)^{{\rm M}(\sigma)} u_3^{{\rm rda}(\sigma)}u_4^{{\rm ldd}(\sigma)}{\beta}^{{\rm RLmin}(\sigma)-1}\nonumber
\end{align}
Given $\sigma=(1,\sigma_2,\ldots, \sigma_{n+1}) \in \mathfrak{S}_{n+1}$. Define $\overline{\sigma}=(\sigma_2-1,\sigma_3-1,\ldots, \sigma_{n+1}-1)$. It is easy to check that $\overline{\sigma} \in \mathfrak{S}_{n}$ and
\[{\rm M}(\sigma)={\rm L}(\overline{\sigma}),\  {\rm rda}(\sigma)={\rm lrda}(\overline{\sigma}),\  {\rm ldd}(\sigma)={\rm dd}(\overline{\sigma}),\ {\rm RLmin}(\sigma)-1={\rm RLmin}(\overline{\sigma}).\]
Moreover, this process is reversible. Hence we derive that
\begin{align}\label{rda-lrdamin}
P_n(u_1,u_2,u_3,u_4|0,{\beta}) =\sum_{{\overline{\sigma}  \in \mathfrak{S}_{n}}}(u_1u_2)^{{\rm L}(\overline{\sigma})} u_3^{{\rm lrda}(\overline{\sigma})}u_4^{{\rm dd}(\overline{\sigma})}{\beta}^{{\rm RLmin}(\overline{\sigma})}.
\end{align}
Using the same argument, we have
\begin{align}\label{defi-Euler-min-0}
A_n(x,y|\alpha,0)&=\sum_{{\sigma} \in \mathfrak{S}_{n}}x^{{\rm des}({\sigma})+1}y^{{\rm asc}({\sigma})}{\alpha}^{{\rm LRmin}({\sigma})}
\end{align}
and
\begin{align}\label{defi-Euler-min-1}
A_n(x,y|0,\beta)&=\sum_{{\sigma} \in \mathfrak{S}_{n}}x^{{\rm des}({\sigma})}y^{{\rm asc}({\sigma})+1}{\beta}^{{\rm RLmin}({\sigma})}
\end{align}

We are now in a position to show  Theorem \ref{thm-grammar-P} by  using the grammatical labeling. The notion of a grammatical labeling was introduced by Chen and Fu \cite{Chen-Fu-2017}.

Let $\sigma=\sigma_1\sigma_2\cdots \sigma_n \in \mathfrak{S}_n$. For $1\leq i\leq n+1$,  recall that the position $i$ is said to be the position immediately before $\sigma_i$, whereas the position $n+1$ is meant to be the position after $\sigma_n$.
The labeling for $P_n(u_1,u_2, u_3,u_4|{\alpha},{\beta})$ can be described as follows. We patch   $+\infty$ to $\sigma$ at both ends so that there are $n+1$ positions between two adjacent elements. For $1\leq i\leq n+1$, we label the position $i$ as follows:
\begin{itemize}
\item If $i=1$, then label it by $a$;
\item  If $i=n+1$, then label it by $b$;
 \item If $i$ is a right double ascent, then label the position $i$  by $u_3$;
 \item If $i-1$ is a left double descent, then label the position $i$  by $u_4$;
 \item If $i$ is a peak, then label the position $i$  by $u_2$ and label the position $i+1$  by $u_1$;
 \item If $\sigma_j$ is a left-to-right minimum and $\sigma_j\neq 1$, then label ${\alpha}$ below  $\sigma_j$;
 \item  If $\sigma_j$ is a right-to-left minimum and $\sigma_j\neq 1$, then label ${\beta}$ below  $\sigma_j$.
\end{itemize}
The weight $\omega$ of $\sigma$ is defined to be the product of all the labels. Below is an example,
\begin{align*}
  & 7\quad 1\quad 2\quad 8\quad 3\quad 6\quad 5\quad 4\\[5pt]
  &{\longrightarrow} \ \begin{array}{ccccccccccccccccccc}
   +\infty&a&7&u_4&1&u_3&2&u_2&8&u_1&3&u_2&6&u_1&5&u_4&4&b&+\infty  \\
     & &{\alpha}&&&&{\beta}&&&&{\beta}&&&&&&{\beta}&&
\end{array}
\end{align*}

If $\sigma$ has $k$  peaks, $d_1$ right double ascents, $d_2$ left double descents,  and $m_1$ left-to-right minima and $m_2$ right-to-left minima, then its weight is given by
\[\omega(\sigma)=ab(u_1u_2)^{k}u_3^{d_1}u_4^{d_2}{\alpha}^{m_1-1}{\beta}^{m_2-1}.\]
From the definition of the above labeling, we see that
\begin{equation}
  ab P_n(u_1,u_2, u_3,u_4|{\alpha},{\beta})=\sum_{\sigma \in \mathfrak{S}_{n+1}}\omega(\sigma).
\end{equation}

\noindent{\bf Proof of Theorem   \ref{thm-grammar-P}.} We proceed by induction on $n$. For $n=0$, the statement is obvious. Assume that this assertion holds for $n-1$. To show that it is valid for $n$, we represent a permutation $\sigma=\sigma_1\sigma_2\cdots \sigma_n$ in $\mathfrak{S}_n$ by patching   $+\infty$  at both ends so that there are $n+1$ positions between two adjacent elements, and these are $n+1$ positions to insert the element $n+1$ into $\sigma$ to generate a permutation in $\mathfrak{S}_{n+1}$. Suppose that $\pi$ is permutation  in $\mathfrak{S}_{n+1}$ obtained from $\sigma$ by inserting the element $n+1$ at the position $i$, where $1\leq i\leq n+1$. We consider the following six cases:

Case 1: If $i$ is labeled by $a$ in the labeling of $\sigma$, that means that $i=1$,  then the first position and the second position of $\pi$ are labeled by $a$ and $u_4$ in the labeling of $\pi$ respectively. Moreover, $n+1$ is the left-to-right minimum of $\pi$, so label ${\alpha}$ below $n+1$. In this case, we find that  the change of weights of $\pi$ is consistent with the substitution rule $a\rightarrow a\alpha u_4$.
\[\begin{array}{cccccc}
+\infty & a&\sigma_1&   \cdots \\
 &  &   {\alpha}&
\end{array}\  \Rightarrow \  \begin{array}{cccccc}
+\infty & a&n+1&u_4&\sigma_1&   \cdots \\
  &  &{\alpha}&&{\alpha}&   \cdots \\
\end{array}.\]

Case 2: If $i$ is labeled by $b$ in the labeling of $\sigma$, that means $i=n+1$, then   the changes of weights caused by the insertion are coded by the rule $b\rightarrow b\beta u_3$.
\[\begin{array}{cccccc}
\cdots &\sigma_n&   b &  +\infty\\
\cdots&      {\beta} &   &
\end{array}\  \Rightarrow \  \begin{array}{cccccc}
\cdots &\sigma_n&  u_3& n+1 & b &  +\infty\\
\cdots&   {\beta} & & {\beta} &   &
\end{array}.\]

Case 3: If $i$ is labeled by $u_3$ in the labeling of $\sigma$, then the position $i$ of $\pi$ is labeled by $u_2$ and the position $i+1$ of $\pi$ is labeled by $u_1$ in the labeling of $\pi$, so the change of weights of $\pi$ is consistent with the substitution rule $u_3\rightarrow u_1u_2$.
\[\cdots\ \sigma_{i-1} \ u_3\quad\sigma_i \  \cdots\  \Rightarrow \cdots  \sigma_{i-1}\ u_2\quad  n+1 \quad u_1\quad  \sigma_i \  \cdots.\]

Case 4: If  $i$ is labeled by $u_4$ in the labeling of $\sigma$, then it can be checked that the change of weights is in accordance with the rule $u_4\rightarrow u_1u_2$, see the figure below.
\[\cdots\  \sigma_{i-1}\quad u_4\quad\sigma_i  \  \cdots\  \Rightarrow\  \cdots \ \sigma_{i-1}\quad u_2 \quad n+1 \quad u_1 \quad  \quad \sigma_{i}  \  \cdots.\]

Case 5: If  $i$ is labeled by $u_1$ in the labeling of $\sigma$, then it can be checked that the change of weights is in accordance with the rule $u_1\rightarrow u_1u_3$, see the figure below.
\[\cdots\    \sigma_{i-2}\quad u_2\quad\sigma_{i-1} \quad u_1 \quad \sigma_i \  \cdots\  \Rightarrow \  \cdots    \sigma_{i-2}\quad u_3\quad\sigma_{i-1} \quad u_2 \quad n+1\quad u_1 \quad \sigma_i \  \cdots\  .\]

Case 6: If  $i$ is labeled by $u_2$ in the labeling of $\sigma$, then it can be checked that the change of weights is in accordance with the rule $u_2\rightarrow u_2u_4$, see the figure below.
\[\cdots\    \sigma_{i-1} \quad u_2 \quad \sigma_i \quad  u_1 \quad  \sigma_{i+1} \  \cdots\  \Rightarrow \  \cdots\    \sigma_{i-1} \quad u_2\quad n+1 \quad u_1 \quad \sigma_i \quad  u_4 \quad  \sigma_{i+1} \  \cdots\ .\]

Summing up all the cases shows that this assertion is valid for $n$. This completes the proof. \qed

\subsection{Proof of  Theorem \ref{thm-grammar-P2}}

We are now in a position to give a grammatical derivation of  Theorem \ref{main2}. By employing Theorem \ref{thm-grammar-P}, it is sufficient to demonstrate Theorem \ref{thm-grammar-P2}.

\noindent {\bf Proof of  Theorem \ref{thm-grammar-P2}.}  Let $D_{\widetilde{G}_1}$ is the formal derivative with respect to the   grammar
 \begin{equation} \label{L-P-four-origi}
    \widetilde{G}_1=\{u_4\rightarrow u_1u_2, \, u_3\rightarrow u_1u_2,\,u_1\rightarrow u_1u_3, u_2\rightarrow  u_2u_4 \}.
     \end{equation}
We first show that
\begin{align}\label{main2-pf-a1}
{\rm Gen}^{(\widetilde{G}_1)}(u_1u_2,t)&:=\sum_{n\geq 0}D_{\widetilde{G}_1}^n(u_1u_2)\frac{t^n}{n!}\nonumber\\[5pt]
&=xy{(1+yF(x,y;t))(1+xF(x,y;t))}
\end{align}
and
\begin{align}\label{main2-pf-a2}
{\rm Gen}^{(\widetilde{G}_1)}(u_1,t)&:=\sum_{n\geq 0}D_{\widetilde{G}_1}^n(u_1)\frac{t^n}{n!}\nonumber\\[5pt]
&=u_1\sqrt{(1+yF(x,y;t))(1+xF(x,y;t))}e^{\frac{(u_3-u_4)t}{2}},
\end{align}
where $x+y=u_3+u_4$ and $xy=u_1u_2$ and $F(x,y;t)$ is given by \eqref{defi-F}.

Recall that $D_{{G}_1}$ is the formal derivative with respect to the   grammar \eqref{C-S-G-orign} and $D_{\widetilde{G}_1} $ is  the formal derivative with respect to the   grammar \eqref{L-P-four-origi}. If we  set  $x+y=u_3+u_4$ and $xy=u_1u_2$,  we find that
\begin{equation}\label{pf-joi-tem-3}
D_{{G}_1}(xy)=xy(x+y)=u_1u_2(u_3+u_4)=D_{\widetilde{G}_1}(u_1u_2)
\end{equation}
and
\begin{equation}\label{pf-joi-tem-4}
D_{{G}_1}(x+y)=2xy=2u_1u_2=D_{\widetilde{G}_1}(u_3+u_4).
\end{equation}
We claim that for $n\geq 1$,
\begin{equation}\label{pf-4variable-rela}
D^n_{{G}_1}(xy)=D^n_{\widetilde{G}_1}(u_1u_2).
\end{equation}
By \eqref{pf-joi-tem-3}, we see that  \eqref{pf-4variable-rela} is valid when $n=1$.  Assume that \eqref{pf-4variable-rela} holds for $n$. Observe that $D^n_{{G}_1}(xy)$ is symmetric in $x,y$, so we may write $D^n_{{G}_1}(xy)$ in the following form:
\begin{equation}\label{pf-joi-tem-1}
D^n_{{G}_1}(xy)=\sum_{j=1}^{\lfloor \frac{n+2}{2}\rfloor}a_{j} (xy)^j (x+y)^{n+2-2j},
\end{equation}
where $a_j$ are positive integers. By the induction hypothesis, we see that
\[D^n_{\widetilde{G}_1}(u_1u_2)=D^n_{{G}_1}(xy). \]
and so
\begin{equation} \label{pf-joi-tem-2}
D^n_{\widetilde{G}_1}(u_1u_2)=\sum_{j=1}^{\lfloor \frac{n+2}{2}\rfloor}a_{j} (u_1u_2)^j (u_3+u_4)^{n+2-2j}
\end{equation}
Applying $D_{{G}_1}$ to \eqref{pf-joi-tem-1}, we obtain that
\begin{align*}
D^{n+1}_{{G}_1}(xy)&=\sum_{j=1}^{\lfloor \frac{n+2}{2}\rfloor}a_{j} j(xy)^{j-1}(x+y)^{n-2j+2} D_{{G}_1}(xy) \nonumber \\[5pt]
&\quad + \sum_{j=1}^{\lfloor \frac{n+2}{2}\rfloor}a_{j} (n-2j+2)(xy)^{j}(x+y)^{n-2j+1} D_{{G}_1}(x+y)
\end{align*}
Since   $x+y=u_3+u_4$ and $xy=u_1u_2$,  we find that
\begin{align}
D^{n+1}_{{G}_1}(xy)
&=\sum_{j=1}^{\lfloor \frac{n+2}{2}\rfloor}a_{j} j(u_1u_2)^{j-1}(u_3+u_4)^{n-2j+2} D_{{G}_1}(xy) \nonumber\\[5pt]
& \quad + \sum_{j=1}^{\lfloor \frac{n+2}{2}\rfloor}a_{j} (n-2j+2)(u_1u_2)^{j}(u_3+u_4)^{n-2j+1} D_{{G}_1}(x+y) \label{pf-joi-tem-5}
\end{align}
Applying \eqref{pf-joi-tem-3} and \eqref{pf-joi-tem-4} into \eqref{pf-joi-tem-5}, and by \eqref{pf-joi-tem-2}, we  conclude that
\[D^{n+1}_{{G}_1}(xy)=D^{n+1}_{\widetilde{G}_1}(u_1u_2),\]
and so \eqref{pf-4variable-rela} also is valid for $n+1$, and hence the claim is verified. Therefore, we obtain
\begin{equation}\label{pf-joi-tem-6}
{\rm Gen}^{(\widetilde{G}_1)}(u_1u_2,t)={\rm Gen}^{({G}_1)}(xy,t).
\end{equation}
By applying the multiplicative property \eqref{gramma-multiple}, we can deduce from \eqref{chen-fu-eulera} and \eqref{chen-fu-eulerb} that
\begin{equation}\label{chen-fu-xy}
{\rm Gen}^{(G_{1})}(xy,t)=\sum_{n\geq 0}D_{{G}_1}^n(xy)\frac{t^n}{n!}=xy{(1+yF(x,y;t))(1+xF(x,y;t))}.
\end{equation}
Substituting \eqref{chen-fu-xy} into \eqref{pf-joi-tem-6}, we obtain \eqref{main2-pf-a1}.

To prove  \eqref{main2-pf-a2}, we first observe that
\begin{equation}
D_{\widetilde{G}_1}(u_1u_2^{-1})=u_1u_2^{-1}(u_3-u_4).
\end{equation}
Since
\[D_{\widetilde{G}_1}(u_3-u_4)=0,\]
it follows that for $n\geq 0$,
\begin{equation*}
D^n_{\widetilde{G}_1}(u_1u_2^{-1})=u_1u_2^{-1}(u_3-u_4)^n.
\end{equation*}
Hence
\begin{align}\label{main2-pf-a222}
{\rm Gen}^{(\widetilde{G}_1)}(u_1u_2^{-1},t)&:=\sum_{n\geq 0}D_{\widetilde{G}_1}^n(u_1u_2^{-1})\frac{t^n}{n!}=u_1u_2^{-1}e^{(u_3-u_4)t}.
\end{align}
By the multiplicative property \eqref{gramma-multiple}, we deduce from \eqref{main2-pf-a1} and \eqref{main2-pf-a222} that
\begin{align*}
({\rm Gen}^{(\widetilde{G}_1)}(u_1,t))^2&={\rm Gen}^{(\widetilde{G}_1)}(u_1u_2^{-1},t){\rm Gen}^{(\widetilde{G}_1)}(u_1u_2,t)\\[5pt]
&=u^2_1((1+yF(x,y;t))(1+xF(x,y;t)))e^{{(u_3-u_4)t}},
\end{align*}
from which, we obtain  \eqref{main2-pf-a2}.

Let ${\alpha}, {\beta}$ be two fixed numbers. We see that
\begin{align}\label{main2-pf-t1}
D_{\widetilde{G}_1}(u_2^{\alpha})&={\alpha}u_2^{{\alpha}-1}
D_{\widetilde{G}_1}(u_2)={\alpha}
u_2^{{\alpha}}u_4, \\[5pt]
D_{\widetilde{G}_1}(u_1^{\beta})&={\beta} u_1^{{\beta}-1}D_{\widetilde{G}_1}(u_1)={\beta}u_1^{{\beta}}u_3.  \label{main2-pf-t2}
\end{align}
Let $D_{\widetilde{G}}$ is the formal derivative with respect to the   grammar \eqref{L-P-four}. Setting $a=u_2^{\alpha}$ and $b=u_1^{\beta}$, then by \eqref{main2-pf-t1} and \eqref{main2-pf-t2}, we find that
\[D_{\widetilde{G}_1}(u_2^{\alpha})=D_{\widetilde{G}}(a) \quad \text{and} \quad
D_{\widetilde{G}_1}(u_1^{\beta})=D_{\widetilde{G}}(b).\]
Moreover, it is easy to check that
\[D_{\widetilde{G}_1}(u_1)=D_{\widetilde{G}}(u_1),  \quad D_{\widetilde{G}_1}(u_2)=D_{\widetilde{G}}(u_2), \quad
D_{\widetilde{G}_1}(u_3)=D_{\widetilde{G}}(u_3), \quad \text{and} \quad D_{\widetilde{G}_1}(u_4)=D_{\widetilde{G}}(u_4).\]
Hence we can use the induction on $n$ to deduce that for $n\geq 0$,
\[D^n_{\widetilde{G}_1}(u_2^{\alpha})=D^n_{\widetilde{G}}(a) \quad \text{and} \quad
D^n_{\widetilde{G}_1}(u_1^{\beta})=D^n_{\widetilde{G}}(b).\]
Consequently,  for $n\geq 0$,
\begin{equation}\label{pf-main2-gf-t3}
D^n_{\widetilde{G}}(ab)=D^n_{\widetilde{G}_1}(u_2^{\alpha}u_1^{\beta}).
\end{equation}
It follows that
\begin{align*}
  {\rm Gen}^{(\widetilde{G})}(ab,t)&= {\rm Gen}^{(\widetilde{G}_1)}(u_2^{\alpha}u_1^{\beta},t)\\[5pt]
&=({\rm Gen}^{({\widetilde{G}_1})}(u_2u_1,t))^{\alpha}({\rm Gen}^{({\widetilde{G}_1})}(u_1,t))^{\beta-\alpha}\\[5pt]
&=u^\beta_1u^\alpha_2\left(1+yF(x,y;t)\right)^\frac{{\alpha}+{\beta}}{2}
\left(1+xF(x,y;t)\right)^\frac{{\alpha}+{\beta}}{2}
e^{\frac{1}{2}({\beta}-{\alpha})(u_3-u_4)t}\\[5pt]
&=ab\left(1+yF(x,y;t)\right)^\frac{{\alpha}+{\beta}}{2}
\left(1+xF(x,y;t)\right)^\frac{{\alpha}+{\beta}}{2}
e^{\frac{1}{2}({\beta}-{\alpha})(u_3-u_4)t}
\end{align*}
as desired. This completes the proof of Theorem \ref{thm-grammar-P2}. \qed

It should be noted that the generating functions \eqref{main2-pf-a1} and \eqref{main2-pf-a2} in the proof of Theorem \ref{main2} could also be derived by using the results of Fu \cite{Fu-2018}.

 \section{The generating functions}

In this section, we give a proof of Theorem \ref{main2c}  with the aid of  Theorem  \ref{main2}. We then derive many consequences of Theorem \ref{main2c}, which provide $(\alpha,\beta)$-extensions of the   generating functions of peaks, left peaks, double  ascents, right double  ascents and left-right double ascents.

\noindent {\bf Proof of  Theorem \ref{main2c}.}  From Theorem \ref{main2}, we see that
 \begin{align}\label{pf2c-tem1}
 &\sum_{n\geq 0}P_n(u_1,u_2,u_3,u_4|{\alpha},{\beta})\frac{t^n}{n!}\nonumber\\[5pt]
 &\quad =\left((1+xF(x,y;t))(1+yF(x,y;t))\right)^\frac{{\alpha}+{\beta}}
 {2}e^{\frac{1}{2}({\beta}-{\alpha})(u_3-u_4)t},
 \end{align}
 where $x+y=u_3+u_4$, $xy=u_1u_2$.

 Recall that
 \[F(x,y;t)=\frac{e^{xt}-e^{yt}}{xe^{yt}-ye^{xt}},\]
we find that
\begin{align}\label{pf2c-tem2}
1+xF(x,y;t)&=\frac{(x-y)e^{xt}}{xe^{yt}-ye^{xt}}\nonumber\\[5pt]
&=\left(-\frac{x}{y-x}e^{(y-x)t}+\frac{y}{y-x}\right)^{-1}.
\end{align}
Similarly, we have
\begin{align}\label{pf2c-tem3}
1+yF(x,y;t)&=\left(\frac{y}{y-x}e^{(x-y)t}
-\frac{x}{y-x}\right)^{-1}.
\end{align}
Since $x+y=u_3+u_4$ and $xy=u_1u_2$, it follows that
\begin{align}
x&=\frac{(u_3+u_4)-\sqrt{(u_3+u_4)^2-4u_1u_2}}{2},\label{id-x-u3u4}\\[5pt]
y&=\frac{(u_3+u_4)+\sqrt{(u_3+u_4)^2-4u_1u_2}}{2}.\label{id-y-u3u4}
\end{align}
Hence, we have
\begin{align}\label{pf2c-tem4}
y-x&=\sqrt{(u_3+u_4)^2-4u_1u_2},
\\[5pt]
\frac{x}{y-x}&=-\frac{1}{2}+\frac{1}{2}
\frac{(u_3+u_4)}{\sqrt{(u_3+u_4)^2-4u_1u_2}}, \label{pf2c-tem5}\\[5pt]
\frac{y}{y-x}&=\frac{1}{2}+
\frac{1}{2}\frac{(u_3+u_4)}{\sqrt{(u_3+u_4)^2-4u_1u_2}}. \label{pf2c-tem6}
\end{align}
Putting \eqref{pf2c-tem4}, \eqref{pf2c-tem5} and \eqref{pf2c-tem6} into \eqref{pf2c-tem2}, we obtain
\begin{align*}
&-\frac{x}{y-x}e^{(y-x)t}+\frac{y}{y-x}\\[5pt]
&\quad =\frac{e^{t\sqrt{(u_3+u_4)^2-4u_1u_2}}+1}{2}-\frac{(u_3+u_4)}{\sqrt{(u_3+u_4)^2-4u_1u_2}}
\frac{e^{t\sqrt{(u_3+u_4)^2-4u_1u_2}}-1}{2} \\[5pt]
&\quad=\frac{1}{e^{-\frac{t}{2}\sqrt{(u_3+u_4)^2-4u_1u_2}}}
\left(\cosh\left(\frac{t}{2}\sqrt{(u_3+u_4)^2-4u_1u_2}\right)\right.\\[5pt]
&\quad\quad \quad\left.
-\frac{(u_3+u_4)}{\sqrt{(u_3+u_4)^2-4u_1u_2}}
\sinh\left(\frac{t}{2}\sqrt{(u_3+u_4)^2-4u_1u_2}\right)\right).
\end{align*}
Similarly, plugging  \eqref{pf2c-tem4}, \eqref{pf2c-tem5} and \eqref{pf2c-tem6} into \eqref{pf2c-tem3}, we derive that
 \begin{align*}
&-\frac{x}{y-x}+\frac{y}{y-x}e^{(x-y)t}\\[5pt]
&\quad=\frac{1}{e^{\frac{t}{2}\sqrt{(u_3+u_4)^2-4u_1u_2}}}
\left(\cosh\left(\frac{t}{2}\sqrt{(u_3+u_4)^2-4u_1u_2}\right)\right.\\[5pt]
&\quad\quad \quad\left.
-\frac{(u_3+u_4)}{\sqrt{(u_3+u_4)^2-4u_1u_2}}
\sinh\left(\frac{t}{2}\sqrt{(u_3+u_4)^2-4u_1u_2}\right)\right).
\end{align*}
Consequently,
 \begin{align}\label{pf2c-tem7}
&(1+xF(x,y;t))(1+yF(x,y;t))=\left(\cosh\left(\frac{t}{2}
\sqrt{(u_3+u_4)^2-4u_1u_2}\right)\right.\nonumber \\[5pt]
&\quad\quad \quad\left.
-\frac{(u_3+u_4)}{\sqrt{(u_3+u_4)^2-4u_1u_2}}
\sinh\left(\frac{t}{2}\sqrt{(u_3+u_4)^2-4u_1u_2}
\right)\right)^{-2}.
\end{align}
Substituting \eqref{pf2c-tem7} into \eqref{pf2c-tem1} yields the generating function of $P_n(u_1,u_2,u_3,u_4|{\alpha},{\beta})$ as stated in Theorem \ref{main2c}. This completes the proof. \qed

Setting $\alpha=0$ in Theorem \ref{main2c}, and by \eqref{rda-lrdamin}
, we have
\begin{thm} \label{cor-lpk-des} We have
\begin{align*}
&\sum_{n\geq 0}\left(\sum_{{\sigma \in \mathfrak{S}_{n}}}(u_1u_2)^{{\rm L}(\sigma)} u_3^{{\rm lrda}(\sigma)}u_4^{{\rm dd}(\sigma)}{\beta}^{{\rm RLmin}(\sigma)}\right)\frac{t^n}{n!}=e^{\frac{\beta}{2}(u_3-u_4) t}\times \left(\cosh\left(\frac{t}{2}\sqrt{(u_3+u_4)^2  -4u_1u_2}\right) \right.\\[5pt]
& \hskip 2.5cm  \left. -\frac{u_3+u_4}{\sqrt{(u_3+u_4)^2-4u_1u_2}}\sinh\left(\frac{t}{2}\sqrt{(u_3+u_4)^2-4u_1u_2}\right)\right)^{-{\beta}}.
\end{align*}
\end{thm}

Setting $u_1=u_2=u$ and $u_3=u_4=v$ in Theorem \ref{cor-lpk-des}, and using \eqref{rel-deslpkdd-a},  we derive the following $\beta$-extension of the   generating function for the left peak polynomials.
\begin{thm}\label{gf-pk-thm} We have
\begin{align}
 &\sum_{n\geq 0}\left( \sum_{{\sigma  \in \mathfrak{S}_{n}}}u^{ 2{\rm L}(\sigma)}v^{n-2{\rm L}(\sigma)} {\beta}^{{\rm RLmin}(\sigma)}\right) \frac{t^n}{n!}\nonumber\\[5pt]
&\quad = \left(\frac{\sqrt{v^2-u^2}}{\sqrt{v^2-u^2} \cosh\left(t\sqrt{v^2-u^2}\right)-v\sinh\left(t\sqrt{v^2-u^2}\right)}\right)^{{\beta}}. \nonumber
\end{align}
\end{thm}
Setting $\beta=1$ in Theorem \ref{gf-pk-thm}, we recover  the   generating function \eqref{gf-lpk} for the left peak polynomials established  by Gessel \cite[Sequence A008971]{OLEN}.

Setting $u_1=u_2=u_4=1$ and $u_3=u$ in Theorem \ref{cor-lpk-des},  we acquire the $\beta$-extension of the    generating function for  left-right double ascents, from which  we recover the    generating function \eqref{zhuang-lrda} for the left-right double ascents by setting  $\beta=1$.
\begin{thm}\label{gf-lrda} We have
\begin{align*}
&\sum_{n\geq 0}\left( \sum_{{\sigma \in \mathfrak{S}_{n}}}u^{ {\rm lrda}(\sigma)}   {\beta}^{{\rm RLmin}(\sigma)}\right) \frac{t^n}{n!}\\[5pt]
&\quad =e^{\frac{\beta (u-1)}{2}t} \left(\frac{\sqrt{(u+1)^2-4}}{\sqrt{(u+1)^2-4} \cosh\left(\frac{t}{2}\sqrt{(u+1)^2-4}\right)-(1+u)\sinh\left(\frac{t}{2}\sqrt{(u+1)^2-4}\right)}\right)^{{\beta}}.
\end{align*}
\end{thm}

Setting $u_1=u_2=u_3=1$ and $u_4=u$ in Theorem \ref{cor-lpk-des}, and by taking reverse of a permutation, we obtain the $\alpha$-extension of the generating function for  double ascents.  This allows us to retrieve the generating function \eqref{Zhuangda} for   double ascents when $\alpha=1$.

\begin{thm}\label{gf-da} We have
\begin{align*}
&\sum_{n\geq 0}\left( \sum_{{\sigma \in \mathfrak{S}_{n}}}u^{ {\rm da}(\sigma)}   {\alpha}^{{\rm LRmin}(\sigma)}\right) \frac{t^n}{n!}\\[5pt]
&\quad =e^{\frac{\alpha(1-u)}{2}t} \left(\frac{\sqrt{(u+1)^2-4}}{\sqrt{(u+1)^2-4} \cosh\left(\frac{t}{2}\sqrt{(u+1)^2-4}\right)-(1+u)\sinh\left(\frac{t}{2}\sqrt{(u+1)^2-4}\right)}\right)^{{\alpha}}.
\end{align*}
\end{thm}

Setting $u_1=u_2=u$ and $u_3=u_4=v$ in Theorem \ref{main2c}  and using \eqref{defi-P-min},
we arrive at

\begin{thm}\label{gf-pk-lr} We have
\begin{align*}
&\sum_{n\geq 0}\left( \sum_{{\sigma \in \mathfrak{S}_{n+1}}}u^{ 2{\rm M}(\sigma)}v^{n-2{\rm M}(\sigma)} {\alpha}^{{\rm LRmin}(\sigma)-1} {\beta}^{{\rm RLmin}(\sigma)-1}\right) \frac{t^n}{n!}\\[5pt]
&\quad = \left(\frac{\sqrt{v^2-u^2}}{\sqrt{v^2-u^2} \cosh\left(t\sqrt{v^2-u^2}\right)-v\sinh\left(t\sqrt{v^2-u^2}\right)}\right)^{{\alpha+\beta}}.
\end{align*}
\end{thm}
By setting $v=\alpha=\beta=1$ in Theorem \ref{gf-pk-lr}, replacing $u$ with $\sqrt{u}$, and performing integration on both sides with respect to $t$, we retrieve  the generating function \eqref{gf-pk} for the peak polynomials.

Setting $u_1=u_2=u_4=1$ and $u_3=u$ in Theorem \ref{main2c} and using \eqref{defi-P-min} gives that

\begin{thm}\label{gf-rda} We have
\begin{align*}
&\sum_{n\geq 0}\left( \sum_{{\sigma \in \mathfrak{S}_{n+1}}}u^{ {\rm rda}(\sigma)}  {\alpha}^{{\rm LRmin}(\sigma)-1} {\beta}^{{\rm RLmin}(\sigma)-1}\right) \frac{t^n}{n!} =e^{\frac{(\beta-\alpha)(u-1)}{2}t}\\[5pt] &\quad \quad \times \left(\frac{\sqrt{(u+1)^2-4}}{\sqrt{(u+1)^2-4} \cosh\left(\frac{t}{2}\sqrt{(u+1)^2-4}\right)-(1+u)\sinh\left(\frac{t}{2}\sqrt{(u+1)^2-4}\right)}\right)^{{\alpha+\beta}}.
\end{align*}
\end{thm}
Setting $\alpha=\beta=1$ in Theorem \ref{gf-rda} and then integrating both sides of the mentioned identity with respect to $t$,  we   bring back the   generating function \eqref{zhuang-rda} for the right double ascents.

\section{The expression of $P_n(u_1,u_2,u_3,u_4|{\alpha},{\beta})$ with $\alpha+\beta=-1$ }

This section is dedicated to deriving  an explicit expression of $P_n(u_1,u_2,u_3,u_4;{\alpha},{\beta})$ when ${\alpha}+{\beta}=-1$ by utilizing  Theorem \ref{main2c}.

\begin{thm} \label{main2ca} When ${\alpha}+{\beta}=-1$ and for $n\geq 1$,
\begin{align} \label{main2caeq}
 & 2^n P_n(u_1,u_2,u_3,u_4|{\alpha},{\beta})  \\[5pt]
  &\quad = \sum_{k=0}^{\lfloor \frac{n}{2}\rfloor} {n\choose 2k} ((u_3+u_4)^2-4u_1u_2)^k({\beta}-{\alpha})^{n-2k}(u_3-u_4)^{n-2k}\nonumber\\[5pt]
& \quad  \quad  -(u_3+u_4)\sum_{k=0}^{\lfloor \frac{n}{2}\rfloor} {n\choose 2k+1} ((u_3+u_4)^2-4u_1u_2)^k({\beta}-{\alpha})^{n-2k-1}(u_3-u_4)^{n-2k-1}. \nonumber
\end{align}
\end{thm}

\pf  When $\alpha+\beta=-1$, and if we let $v=\sqrt{(u_3+u_4)^2-4u_1u_2}$,  then Theorem \ref{main2c} becomes
\begin{align*}
&\sum_{n\geq 0} P_n(u_1,u_2,u_3,u_4|{\alpha},{\beta})  \frac{t^n}{n!}\nonumber\\[5pt]
&=e^{\frac{1}{2}({\beta}-{\alpha})(u_3-u_4)t}\times \left(\cosh\left(\frac{1}{2}vt\right)  -\frac{u_3+u_4}{v}\sinh\left(\frac{1}{2}vt\right)\right)\nonumber\\[5pt]
&=e^{\frac{1}{2}({\beta}-{\alpha})(u_3-u_4)t}\times \left(\sum_{n\geq 0}   \frac{v^{2n}}{2^{2n}}  \frac{t^{2n}}{(2n)!}-(u_3+u_4)\left(\sum_{n\geq 0}  \frac{v^{2n}}{2^{2n+1}}\frac{t^{2n+1}}{(2n+1)!}\right)\right)\nonumber\\[5pt]
&=\left(\sum_{n\geq 0} \frac{(\beta-\alpha)^n(u_3-u_4)^n}{2^n}\frac{t^n}{n!}\right)\times \left(\sum_{n\geq 0}   \frac{v^{2n}}{2^{2n}}  \frac{t^{2n}}{(2n)!}-(u_3+u_4)\left(\sum_{n\geq 0}  \frac{v^{2n}}{2^{2n+1}}\frac{t^{2n+1}}{(2n+1)!}\right)\right)\nonumber\\[5pt]
 &= \sum_{n\geq 0}\frac{t^n}{n!}\left(\sum_{k\geq 0} {n\choose 2k}(\beta-\alpha)^{n-2k}(u_3-u_4)^{n-2k}\frac{v^{2k}}{2^n} \right. \nonumber\\[5pt]
 &\quad \quad \left.-(u_3+u_4) \sum_{k\geq 0}{n\choose 2k+1}(\beta-\alpha)^{n-2k-1}(u_3-u_4)^{n-2k-1}  \frac{v^{2k}}{2^{n}}\right).
\end{align*}
Comparing the coefficients of $t^n/n!$ on the both sides yields \eqref{main2caeq}. This completes the proof. \qed

Setting $u_3=u_4=u_1=1$ and $u_2=u$  in  Theorem \ref{main2ca} and using \eqref{defi-P-min}  yields Theorem \ref{enumer-pk}.  Theorem \ref{enumer-lpk} follows from Theorem \ref{enumer-pk} by setting $\alpha=0$.

By choosing ${\alpha}={\beta}=-1/2$ in Theorem \ref{main2ca}, we derive that
\begin{thm} \label{main2cacor3}For $n\geq 1$,
\begin{align*}
 & 2^n P_n\left(u_1,u_2,u_3,u_4|-\frac{1}{2},-\frac{1}{2}\right) \\[5pt]
  &\quad =\begin{cases} ((u_3+u_4)^2-4u_1u_2)^{\lfloor \frac{n}{2}\rfloor},  \quad \text{if \ $n$ is even},  \\[5pt]
-(u_3+u_4)((u_3+u_4)^2-4u_1u_2)^{\lfloor \frac{n}{2}\rfloor},  \quad \text{if \ $n$ is odd}.
\end{cases}
\end{align*}
\end{thm}
 Setting $u_1=u_2=u_4=1$, $u_3=u$   in  Theorem \ref{main2cacor3} and employing \eqref{defi-P-min} yields Theorem \ref{enumer-rda}.

\section{Relations  between $P_n(u_1,u_2,u_3,u_4|{\alpha},{\beta})$ and $A_n(x,y|{\alpha},{\beta})$ }

In this section, we   begin by presenting two relations between
$P_n(u_1,u_2,u_3,u_4|{\alpha},{\beta})$ and $A_n(x,y|{\alpha},{\beta})$ and give their proofs. Subsequently, we  derive several consequences from these connections. These specific derivations will not only yield   the $(\alpha,\beta)$-extensions of the  related relations associated with the Eulerian polynomial due to Stembridge, Petersen, Br\"and\'en and Zhuang, but will also provide the left peak version of Stembridge's formula,  the peak version of Petersen's formula and their $(\alpha,\beta)$-extensions.

\subsection{Two relations}
\begin{thm} \label{main2-cora} For $n\geq 0$,
\begin{equation}\label{C-S-3-varia-g}
P_n(u_1,u_2,u_3,u_4|{\alpha},{\beta})=\sum^n_{k=0}{n\choose k} A_k\left(x,y|\frac{{\alpha}+{\beta}}{2}, \frac{{\alpha}+{\beta}}{2}\right)\frac{({\beta}-{\alpha})^{n-k}(u_3-u_4)^{n-k}}{2^{n-k}},
\end{equation}
where $x+y=u_3+u_4$ and $xy=u_1u_2$.
\end{thm}
\pf Combining Theorem \ref{C-S-main1} and Theorem \ref{main2}, we derive that
 \begin{align}\label{pf2c-tem1a}
 &\sum_{n\geq 0}P_n(u_1,u_2,u_3,u_4|{\alpha},{\beta})\frac{t^n}{n!}
 \nonumber\\[5pt]
 &\quad =
 e^{\frac{1}{2}({\beta}-{\alpha})(u_3-u_4)t}\sum_{k\geq 0}A_k\left(
 x,y|\frac{{\alpha}+{\beta}}{2},\frac{{\alpha}+{\beta}}{2}
 \right)\frac{t^k}{k!}\nonumber\\[5pt]
 &\quad =\left(\sum_{m\geq 0} \frac{(\beta-\alpha)^m(u_3-u_4)^m}{2^m}\frac{t^m}{m!}\right)\left(
 \sum_{k\geq 0}A_k\left(
 x,y|\frac{{\alpha}+{\beta}}{2},\frac{{\alpha}+{\beta}}{2}
 \right)\frac{t^k}{k!}\right) \nonumber\\[5pt]
  &\quad =\sum_{n\geq 0}\frac{t^n}{n!}\left(\sum^n_{k=0}{n\choose k}A_k\left(
 x,y\Large|\frac{{\alpha}+{\beta}}{2},\frac{{\alpha}+{\beta}}{2}
 \right) \frac{(\beta-\alpha)^{n-k}(u_3-u_4)^{n-k}}{2^{n-k}}\right).
 \end{align}
 Equating the coefficients of $t^n/n!$ yields the result. \qed

 \begin{thm} \label{main2-corb} For $n\geq 0$,
\begin{align}\label{C-S-3-varia-g1}
&P_n(u_1,u_2,u_3,u_4|{\alpha},{\beta})\nonumber\\[5pt]
&\quad =\sum^n_{k=0}{n\choose k} A_k\left(x,y|0, {\alpha}+{\beta}\right)(\alpha x-\beta y+(\beta-\alpha)u_3)^{n-k}\\[5pt]
&\quad =\sum^n_{k=0}{n\choose k} A_k\left(x,y| {\alpha}+{\beta},0\right)(\alpha y-\beta x+(\beta-\alpha)u_3)^{n-k} \label{C-S-3-varia-g2}
\end{align}
where $x+y=u_3+u_4$ and $xy=u_1u_2$.
\end{thm}

\pf  Since, by  \eqref{defi-F}, we see that
  \[1+xF(x,y;t)=
 \frac{(x-y)e^{xt}}{xe^{yt}-ye^{xt}}\]
 and
 \[1+yF(x,y;t)=
 \frac{(x-y)e^{yt}}{xe^{yt}-ye^{xt}}.\]
Hence
 \[(1+xF(x,y;t))(1+yF(x,y;t))=
 \frac{(x-y)^2e^{(x+y)t}}{(xe^{yt}-ye^{xt})^2},\]
so that
 \begin{align}\label{pf-main2-corb1}
&\left((1+xF(x,y;t))(1+yF(x,y;t))\right)^\frac{{\alpha}+{\beta}}
 {2}e^{\frac{1}{2}({\beta}-{\alpha})(u_3-u_4)t}\nonumber\\[5pt]
 &=\left(\frac{x-y}
 {xe^{yt}-ye^{xt}}\right)^{{\alpha}+{\beta}}
 e^{\left(\frac{({\beta}-{\alpha})(u_3-u_4)}{2}
 +\frac{(\alpha+\beta)(x+y)}{2}\right)t}\nonumber\\[5pt]
 &=(1+yF(x,y;t))^{{\alpha}+{\beta}}e^{\left(\frac{({\beta}-{\alpha})(u_3-u_4)}{2}
 +\frac{(\alpha+\beta)(x-y)}{2}\right)t} \\[5pt]
 &=(1+xF(x,y;t))^{{\alpha}+{\beta}}e^{\left(\frac{({\beta}-{\alpha})(u_3-u_4)}{2}
 +\frac{(\alpha+\beta)(y-x)}{2}\right)t} \label{pf-main2-corb2}
 \end{align}
 Applying \eqref{C-S-main1-eq} into \eqref{pf-main2-corb1}, and applying   Theorem \ref{main2},  we find that
 \begin{align*}
 &\sum_{n\geq 0}P_n(u_1,u_2,u_3,u_4|{\alpha},{\beta})\frac{t^n}{n!}
 \nonumber\\[5pt]
 &\quad =(1+yF(x,y;t))^{{\alpha}+{\beta}}e^{\left(\frac{({\beta}-{\alpha})(u_3-u_4)}{2}
 +\frac{(\alpha+\beta)(x-y)}{2}\right)t}\nonumber\\[5pt]
 &\quad =
\left(\sum_{k\geq 0}A_k\left(
 x,y|0,\alpha+\beta
 \right)\frac{t^k}{k!}\right)\left(\sum_{m\geq 0} {\left(\alpha x-\beta y+(\beta-\alpha)u_3\right)}^m \frac{t^m}{m!}\right),
 \end{align*}
 where the last line follows from the relation $x+y=u_3+u_4$.  Equating the coefficients of $t^n/n!$ yields the relation \eqref{C-S-3-varia-g1}.

 If we plug \eqref{C-S-main1-eq} into \eqref{pf-main2-corb1}, and by Theorem \ref{main2},   we derive that
 \begin{align*}
 &\sum_{n\geq 0}P_n(u_1,u_2,u_3,u_4|{\alpha},{\beta})\frac{t^n}{n!}
 \nonumber\\[5pt]
&\quad =(1+xF(x,y;t))^{{\alpha}+{\beta}}e^{\left(\frac{({\beta}-{\alpha})(u_3-u_4)}{2}
 +\frac{(\alpha+\beta)(y-x)}{2}\right)t}\nonumber\\[5pt]
 &\quad =
\left(\sum_{k\geq 0}A_k\left(
 x,y|\alpha+\beta,0
 \right)\frac{t^k}{k!}\right)\left(\sum_{m\geq 0} {\left( \alpha y-\beta x+(\beta-\alpha)u_3 \right)}^m \frac{t^m}{m!}\right).
 \end{align*}
 Equating the coefficients of $t^n/n!$ yields the relation \eqref{C-S-3-varia-g2}. This completes the proof of Theorem \ref{main2-corb}.   \qed

 \subsection{Some consequences}

Setting $u_3=u_4$ and $u_1=u_2$ in  Theorem \ref{main2-cora}, we find that
\begin{equation}\label{relation1}
u_3=u_4=\frac{x+y}{2}\quad \text{and} \quad  u_1=u_2=\sqrt{xy},
\end{equation}
and using the combinatorial definition \eqref{defi-A-min} of  $A_n(x,y|\alpha,\beta)$ and the combinatorial definition \eqref{defi-P-min} of  $P_n(u_1,u_2,u_3,u_4|{\alpha},{\beta})$, we obtain Theorem \ref{main2-cora-1}. Setting $\alpha=\beta=1$ and $y=1$ in Theorem \ref{main2-cora-1},  we reacquire the relation \eqref{Stembridge-v} due to Stembridge.

Choosing $\alpha=0$ in Theorem \ref{main2-cora-1} yields  Theorem \ref{main2-cora-2a}. Setting $\beta=1$ and $y=1$ in Theorem \ref{main2-cora-2a},  we obtain the following consequence, which can be viewed as the left peak version of Stembridge's formula.

\begin{thm}\label{leftpk-Stembridge} For $n\geq 1$,
\begin{equation}\label{Stembridge-v-lpk}
 \sum_{\sigma \in \mathfrak{S}_{n+1}} x^{{\rm des}(\sigma)}{\left(\frac{{1}}{2}\right)}^{{\rm LRmin}(\sigma)+{\rm RLmin}(\sigma)-2}=\left(\frac{1+x}{2}\right)^{n} \sum_{\sigma \in \mathfrak{S}_{n}}\left(\frac{4x}{(1+x)^2}\right)^{{\rm L}(\sigma)}.
 \end{equation}
 \end{thm}

Setting ${\alpha}={\beta}$, $u_1=uv$, $u_2=u_3=w$ and $u_4=v$ in  Theorem \ref{main2-cora}, and by \eqref{defi-A-min}  and \eqref{defi-P-min}, and invoking  \eqref{rel-despkdd-a}, \eqref{rel-compa} \eqref{rel-compb} and \eqref{rel-compd}, we deduce the following consequence, which can be viewed as the $\alpha$-extension of Zhuang's relation \eqref{zhuang-pk-des}.
\begin{thm} \label{zhuang-pk-ex}For $n\geq 1$,
\begin{align} \label{main2-cora-2}
&\sum_{\sigma \in \mathfrak{S}_{n}}u^{{\rm M}(\sigma)}v^{{\rm des}(\sigma)}w^{{\rm asc}(\sigma)}{\alpha}^{{\rm LRmin}(\sigma)+{\rm RLmin}(\sigma)-2}\nonumber \\[5pt]
&=\sum_{\sigma \in \mathfrak{S}_{n}}x^{{\rm des}(\sigma)}y^{{\rm asc}(\sigma)}{{\alpha}}^{{\rm LRmin}(\sigma)+{\rm RLmin}(\sigma)-2},
\end{align}
where
\begin{equation}\label{defi-x}
x=\frac{(w+v)-\sqrt{(w+v)^2-4uvw}}{2},
\end{equation}
and
\begin{equation}\label{defi-y}
y=\frac{(w+v)+\sqrt{(w+v)^2-4uvw}}{2}.
\end{equation}
\end{thm}
Note that \eqref{defi-x} and \eqref{defi-y} follows from \eqref{id-x-u3u4} and \eqref{id-y-u3u4} upon setting $u_1=uv$, $u_2=u_3=w$ and $u_4=v$. Setting ${\alpha}=1$, $w=1$,  $a=x/y$ and $b=(y-1)/(1-x)$ in Theorem \ref{zhuang-pk-ex},  we regain the relation \eqref{zhuang-pk-des} established by  Zhuang \cite[Theorem 4.2]{Zhuang-2017}.

By choosing $\alpha=0$, $u_3=u_4$ and $u_1=u_2$ in \eqref{C-S-3-varia-g2} of Theorem \ref{main2-corb}, and using  \eqref{rda-lrdamin} and \eqref{defi-Euler-min-0},  we derive that the $\beta$-extension of Petersen's relation, from which  we recover the  relation \eqref{relation-lpk}  by setting $y=1$ and $\beta=1$.

 \begin{thm}\label{rela-lpk-exte} For $n\geq 1$,
\begin{align}
&\sum_{\sigma \in \mathfrak{S}_{n}}(xy)^{{\rm L}(\sigma)}\left(\frac{x+y}{2}\right)^{n-2{\rm L}(\sigma)}{\beta}^{{\rm RLmin}(\sigma)}\nonumber \\[5pt]
&=\sum^n_{k=1}{n\choose k} \frac{({\beta}(y-x))^{n-k}}{2^{n-k}}\left(
\sum_{\sigma \in \mathfrak{S}_{k}}x^{{\rm des}(\sigma)+1}y^{{\rm asc}(\sigma)}{\beta}^{{\rm LRmin}(\sigma)}\right)\nonumber \\[5pt]
&\quad \quad \quad + \left(\frac{({\beta}(y-x))}{2}\right)^n.
\end{align}
\end{thm}

  Setting $u_3=u_4$ and $u_1=u_2$ in \eqref{C-S-3-varia-g2} of Theorem \ref{main2-corb}, and using   \eqref{defi-P-min} and \eqref{defi-Euler-min-0},  we get

  \begin{thm}\label{main-coraa}For $n\geq 1$,
\begin{align}
&\sum_{\sigma \in \mathfrak{S}_{n+1}}(xy)^{{\rm M}(\sigma)}\left(\frac{x+y}{2}\right)^{n-2{\rm M}(\sigma)}{\alpha}^{{\rm LRmin}(\sigma)-1}{\beta}^{{\rm RLmin}(\sigma)-1}\nonumber \\[5pt]
&=\sum^n_{k=1}{n\choose k} \frac{(({\alpha}+{\beta})(y-x))^{n-k}}{2^{n-k}}\left(
\sum_{\sigma \in \mathfrak{S}_{k}}x^{{\rm des}(\sigma)+1}y^{{\rm asc}(\sigma)}{(\alpha+\beta)}^{{\rm LRmin}(\sigma)}\right)\nonumber \\[5pt]
&\quad \quad \quad + \left(\frac{(({\alpha}+{\beta})(y-x))}{2}\right)^n.
\end{align}
\end{thm}

By setting $y=1$ and $\alpha=\beta=1$ in Theorem \ref{main-coraa}, we derive the following relation, which can be viewed as    the peak
version of Petersen's formula \eqref{relation-lpk}.
 \begin{thm}\label{pk-Peterson} For $n\geq 1$,
\begin{align*}
&\left(\frac{1+x}{2}\right)^n\sum_{\sigma \in \mathfrak{S}_{n+1}}\left(\frac{4x}{(1+x)^2}\right)^{{\rm M}(\sigma)}\\[5pt]
&\quad \quad =\sum_{k=1}^n {n\choose k}  (1-x)^{n-k}\sum_{\sigma \in \mathfrak{S}_{k}}x^{{\rm des}(\sigma)+1}2^{{\rm LRmin}(\sigma)}+(1-x)^n.
\end{align*}
 \end{thm}

Setting $u_1=uv$, $u_2=u_3=w$ and $u_4=v$ in \eqref{C-S-3-varia-g2} of  Theorem \ref{main2-corb}, and using \eqref{defi-P-min} and   \eqref{defi-Euler-min-0}, we obtain the following consequence. As we will see, this result can be viewed as the $(\alpha,\beta)$-extensions of the peak version of Zhuang's relation \eqref{zhuang-lpk-des}.

\begin{thm}  \label{main2-coraathm}  For $n\geq 1$,
\begin{align} \label{main2-coraa}
&\sum_{\sigma \in \mathfrak{S}_{n+1}}u^{{\rm M}(\sigma)}v^{{\rm des}(\sigma)}w^{{\rm asc}(\sigma)}{\alpha}^{{\rm LRmin}(\sigma)-1}{\beta}^{{\rm RLmin}(\sigma)-1}\\[5pt]
&=\sum^n_{k=1}{n\choose k} \left(\alpha y-\beta x+(\beta-\alpha)w\right)^{n-k}\left(
\sum_{\sigma \in \mathfrak{S}_{k}}x^{{\rm des}(\sigma)+1}y^{{\rm asc}(\sigma)}{(\alpha+\beta)}^{{\rm LRmin}(\sigma)}\right)\nonumber \\[5pt]
&\quad \quad \quad + \left(\alpha y-\beta x+(\beta-\alpha)w\right)^{n},
\end{align}
where $x$ and $y$ are given by \eqref{defi-x} and \eqref{defi-y} respectively.
\end{thm}

The following consequence follows from Theorem \ref{main2-coraathm}  by setting $\alpha=0$.

\begin{thm} \label{rela-zhuang-lpk-exte} For $n\geq 0$,
\begin{align}
&\sum_{\sigma \in \mathfrak{S}_{n}}u^{{\rm L}(\sigma)}v^{{\rm des}(\sigma)}w^{n-{\rm des}(\sigma)} {\beta}^{{\rm RLmin}(\sigma)}\nonumber \\[5pt]
&=\sum^n_{k=1}{n\choose k}\left(\beta(w-x)\right)^{n-k}\left(
\sum_{\sigma \in \mathfrak{S}_{k}}x^{{\rm des}(\sigma)+1}y^{{\rm asc}(\sigma)}{\beta}^{{\rm LRmin}(\sigma)}\right)+ \left(\beta (w-x)\right)^{n},
\end{align}
where $x$ and $y$ are given by \eqref{defi-x} and \eqref{defi-y} respectively.
\end{thm}
Setting ${\beta}=1$, $w=1$,  $a=x/y$ and $b=(y-1)/(1-x)$ in Theorem \ref{rela-zhuang-lpk-exte}, we recover the relation \eqref{zhuang-lpk-des}  due to Zhuang \cite[Theorem 4.7]{Zhuang-2017}.

 \vskip 0.2cm
\noindent{\bf Acknowledgment.} This work
was supported by   the National Science Foundation of China.

\end{document}